\documentclass[11pt, a4paper]{article}
\usepackage{graphicx}
\usepackage{color}
\usepackage{amsmath}
\usepackage{amssymb}
\usepackage{amscd}
\usepackage{bbm}
\usepackage{tikz}
\usepackage{hyperref}
\hypersetup{
    bookmarks=true,         
    colorlinks=true,       
    linkcolor=red,          
    citecolor=green,        
    filecolor=magenta,      
    urlcolor=cyan           
}
\usepackage{marginnote}

\newcommand{\R}{\mathbb{R}}
\newcommand{\inr}[1]{\bigl< #1 \bigr>}

\newcommand{\N}{\mathbb{N}}
\newcommand{\cN}{\mathcal{N}}

\newcommand{\E}{\mathbb{E}}

\newcommand{\eps}{\varepsilon}

\newcommand{\cX}{{\cal X}}
\newcommand{\cH}{{\cal H}}
\newcommand{\cE}{{\cal E}}
\newcommand{\cA}{{\cal A}}

\newcommand{\cO}{{\cal O}}

\newtheorem{Theorem}{Theorem}[section]
\newtheorem{Lemma}[Theorem]{Lemma}

\newtheorem{Definition}[Theorem]{Definition}

\newtheorem{Proposition}[Theorem]{Proposition}
\newtheorem{Corollary}[Theorem]{Corollary}
\newtheorem{Remark}[Theorem]{Remark}

\newtheorem{Assumption}{Assumption}[section]
\newtheorem{Question}[Theorem]{Question}
\numberwithin{equation}{section}

\newcommand{\norm}[1]{\left\|#1\right\|}%

\newcommand{\cD}{{\mathcal{D}}}

\DeclareMathOperator*{\argmin}{argmin}

\def \proof {\noindent {\bf Proof.}\ \ }

\def \endproof
{{\mbox{}\nolinebreak\hfill\rule{2mm}{2mm}\par\medbreak}}

\begin{document}
\title{{Regularization and the small-ball method II: complexity dependent error rates}}
\author{Guillaume Lecu\'e${}^{1,3,5}$  \and Shahar Mendelson${}^{2,4,6}$}

\footnotetext[1]{CNRS, CREST, ENSAE, Bureau E31, 3 avenue Pierre Larousse, 92245 Malakoff.}
\footnotetext[2]{Department of Mathematics, Technion, I.I.T., Haifa, Israel and Mathematical Sciences Institute, The Australian National University, Canberra, Australia}
 \footnotetext[3] {Email: guillaume.lecue@ensae.fr }
\footnotetext[4] {Email: shahar@tx.technion.ac.il}
\footnotetext[5]{Supported by Chaire Havas-Dauphine "Economie des nouvelles donn\'{e}es" and by Investissements d'Avenir (ANR-11-IDEX-0003/Labex Ecodec/ANR-11-LABX-0047)}
\footnotetext[6]{Supported by the Israel Science Foundation, grant 707/14.}

\maketitle

\begin{abstract}
For a convex class of functions $F$, a regularization functions $\Psi(\cdot)$ and given the random data $(X_i, Y_i)_{i=1}^N$, we study estimation properties of regularization procedures of the form
\begin{equation*}
  \hat f \in\argmin_{f\in
    F}\Big(\frac{1}{N}\sum_{i=1}^N\big(Y_i-f(X_i)\big)^2+\lambda \Psi(f)\Big)
\end{equation*}
for some well chosen regularization parameter $\lambda$.

We obtain bounds on the $L_2$ estimation error rate that depend on the complexity of the ``true model" $F^*:=\{f\in F: \Psi(f)\leq\Psi(f^*)\}$, where $f^*\in\argmin_{f\in F}\E(Y-f(X))^2$ and the $(X_i,Y_i)$'s are independent and distributed as $(X,Y)$. Our estimate holds under weak stochastic assumptions  -- one of which being a small-ball condition satisfied by $F$ -- and for  rather flexible choices of regularization functions $\Psi(\cdot)$. Moreover, the result holds in the learning theory framework: we do not assume any a-priori connection between the output $Y$ and the input $X$).

As a proof of concept, we apply our general estimation bound to various choices of $\Psi$, for example, the $\ell_p$ and $S_p$-norms (for $p\geq1$), weak-$\ell_p$, atomic norms, max-norm and SLOPE. In many cases, the estimation rate almost coincides with the minimax rate in the class $F^*$.
\end{abstract}

\section{Introduction}
\label{sec:introduction}
In the standard learning framework, one would like to approximate / predict an unknown random variable $Y$ using functions from a given class $F$, and to do so using only random data. To be more accurate, let $(\cX,\mu)$ be a probability space and consider a class of functions $F$ on $(\cX,\mu)$. Let $X$ be distributed according to $\mu$ and set $X_1,\ldots,X_N \in \cX$ to be $N$ independent copies of $X$.

Given an unknown random variable $Y$, let $\cD=(X_i,Y_i)_{i=1}^N$ be a sample selected according to the joint distribution of $(X,Y)$. One would like to use the data $\cD$ and construct a (random) function $\hat f(\cdot)=\hat f(\cD,\cdot) \in F$, with $\hat f(X)$ serving as a good guess of $Y$.

While there are various interpretations of the meaning of `a good guess', the notion we will focus on here is as follows.

In a typical problem, very little is assumed on the target $Y$ or on the measure $\mu$; on the other hand, the class $F$ is known and a typical assumption is that \textbf{$F$ is convex and closed} in $L_2(\mu)$. Therefore, the functional $f \to \E(f(X)-Y)^2$ has a unique minimizer in $F$,
\begin{equation}\label{eq:oracle}
  f^*=\argmin_{f\in F}\E(Y-f(X))^2.
\end{equation}

The notion of `a good guess' is that $\hat f$ is close to $f^*$ in $L_2(\mu)$, and one would like to obtain a high probability bound on the $L_2(\mu)$ distance of the form \begin{equation}
  \label{eq:risk}
  \norm{\hat f-f^*}_{L_2}^2=\E\left[\big(f^*(X)-\hat f(X)\big)^2| \cD \right] \leq \alpha_N^2.
\end{equation}
In this case,  $\alpha_N^2$ is called a \textit{rate of convergence}, the \textit{error rate} or the \textit{$L_2(\mu)$-estimation rate} of the problem.

\vskip0.4cm

Clearly, one has to pay a price for allowing a rather general target $Y$. Also, to have any hope that $f^*$ is reasonably close to $Y$, one has to consider large classes, leading to an error $\alpha_N^2$ that is often too large to be of any use.

A possible way of bypassing the fact that $F$ may be very large, is the classical approach to {\it regularization}, where a certain property one believes $f^*$ to possess is emphasized by penalizing functions that do not have that property. The penalty is endowed via a {\it regularization function} $\Psi(\cdot)$, defined on an appropriate subspace $E \subset L_2(\mu)$ that contains $F$, and for which $\Psi(f^*)$ is believed to be small (though one does not know that for certain). As a consequence, regularization procedures are designed to fit the data and to have a small $\Psi$ value at the same time. One way of achieving that is to search for functions in $F$ that realize a good trade-off between fitting that data, which is measured via an empirical loss function $P_N\ell_f$, and the size of the regularization term $\lambda \Psi(f)$.

\begin{Definition}\label{def:reg_procedure} The Regularized Empirical Risk Minimization procedure (RERM) is defined by
\begin{equation}
  \label{eq:RERM}
  \hat f \in\argmin_{f\in F}\big(P_N \ell_f+\lambda \Psi(f)\big),
\end{equation}
where here and throughout the article, $P_N h$ denotes the empirical mean of $h$, $\ell_f$ is the loss function associated with $f$ and $\lambda$ is the so-called
regularization parameter.
\end{Definition}

We only consider the square loss $\ell_f(x,y)=(y-f(x))^2$, and thus,
$$
P_N\ell_f=\frac{1}{N}\sum_{i=1}^N(Y_i-f(X_i))^2.
$$

A well known example to this, the ``classical approach" to regularization, is the cubic smoothing spline that can be obtained with the choice 
\begin{equation*}
  \Psi(f)=\int f^{\prime\prime}(t) dt.
\end{equation*}
Another well-studied example is of the form
\begin{equation*}
  \Psi(f)=\int_{\R^d} \frac{\bar{f}(t)}{\bar{G}(t)}dt
\end{equation*}
where the integration is with respect to the Lebesgue measure, $\bar{f}$ is the Fourier transform of $f$ and
$\bar{G}$ is some positive function tending to zero when $|t|$ goes to
infinity (cf. \cite{girosi}). In fact, this type of regularization methods dates back to Tikhonov
(\cite{MR0009685}) and is sometimes called Tikhonov regularization; it is also known as
$L_2$-regularization or Ridge regularization (\cite{MR533250}). 

These methods and others like them have been used to ``smooth'' estimators that have poor generalization capability because of their tendency to over-fit the data, and for the corresponding regularization functions, having a small $\Psi$ value is a guarantee of smoothness.  We refer to \cite{MR2722294} for other examples of regularization functions that have been used to ``smooth'' estimators.

We said ``classical approach to regularization'' because in the more modern approach the aim is somewhat different. One uses a penalty that seemingly has little to do with the property one wishes to emphasize (usually, some notion of {\it sparsity}). Yet somehow, almost ``magically'', the penalty enhances a hidden property and the resulting error rate does not depend on $\Psi(f^*)$ but on that hidden property of $f^*$. We call such error rates \textbf{sparsity-dependent error rates}.

The first part of this article (\cite{LM_sparsity}) has dealt with the modern approach to regularization. Here we would like to complete the picture by exploring bounds that depend on $\Psi(f^*)$ rather than on some hidden sparsity structure of $f^*$. Such error rates will be called \textbf{complexity-dependent error rates}, since the aim is to obtain rates of convergence that depend on the complexity of the unknown ``true model'' $\{f\in F: \Psi(f)\leq \Psi(f^*)\}$. Of course,  the two approaches may sometimes be combined advantageously (see some examples below).

In this context, we will consider regularization functions that satisfy the following properties, which are more general than the ones considered in \cite{LM_sparsity}.

\begin{Assumption} \label{assum:reg-function}
A function $\Psi:E \to \R_{+}$ is a regularization function if
\begin{description}
\item{$\bullet$} It is nonnegative, even, convex and $\Psi(0)=0$.
\item{$\bullet$} There is a constant $\eta \geq 1$, for which, for every $f,h \in E$,
$$
\Psi(f+h) \leq \eta(\Psi(f)+\Psi(h)).
$$
\item{$\bullet$} For every $0 \leq \alpha \leq 1$ and $h \in E$, $\Psi(\alpha h) \leq \alpha \Psi(h) $.
\end{description}
\end{Assumption}

\begin{Remark}
Classical {\it Model Selection regularization
functions}, such as the cardinality of the support of a vector or the rank of
a matrix, are usually not convex and do not satisfy Assumption~\ref{assum:reg-function}. Such examples are therefore not considered in what follows.
\end{Remark}

\subsection{Classical vs. modern}
\label{sec:lasso_intro}
As mentioned above, the direction we take here is closely related to the classical approach to regularization and is rather different from the modern approach. To explain the differences we shall use the celebrated LASSO estimator (cf. \cite{MR1379242,MR1278886}) as an example.

\vskip0.4cm

Let $F$ be a class of linear functionals on $\R^d$ of the form $\inr{t,\cdot}$. Set  $t^*\in\argmin_{t\in \R^d}\E(Y-\inr{X,t})^2$, and consider the RERM \eqref{eq:RERM} with the $\ell_1^d$-norm, $\|t\|_1=\sum_{i=1}^d |t_i|$, serving as a regularization function. Let 
\begin{equation*}
\hat t \in \argmin_{t\in\R^d}\Big(\frac{1}{N}\sum_{i=1}^N \big(Y_i-\inr{X_i,t}\big)^2 + \lambda \norm{t}_1\Big),
\end{equation*}
and the resulting minimizer is the LASSO estimator.

Estimation, de-noising, prediction and support recovery results  have been
obtained for the LASSO in the last decades (see, for example, \cite{MR1379242}, \cite{MR2533469}, and the books \cite{MR3307991, MR2807761} and \cite{MR2829871}
for additional references).
\vskip0.4cm

The LASSO has been used in `high-dimensional' problems, in which the aim was to enhance a low-dimensional structure. The hope was that if the signal $t^*$ were sparse  (that is, supported on relatively few coordinates), the regularization procedure $\hat t$ would estimate $t^*$ with an error rate depending on the cardinality of the support of $t^*$, denoted by $\norm{t^*}_0=|\{j\in\{1,\ldots,d\}: t^*_j\neq0\}|$.

However, if $t^*$ happens to be `well-spread' rather than sparse, though with a reasonable $\ell_1^d$ norm, the sparsity-dependent error rate is useless, while a complexity-dependent error rate, which yields bounds in terms of $\norm{t^*}_1$, is sharper. The obvious example is  $t_1^*=(1/d,...,1/d)$ and $t_2^*=(1,0,...,0)$: although $\|t_1^*\|_1= \|t_2^*\|_1=1$, the cardinalities of their supports are very different, and sparsity-dependent error rates when $t^*=t_1^*$ are likely to be bad.

Examples of that nature are the reason why error rates combining both sparsity and complexity have been obtained for the LASSO. A typical example is Corollary~9.1 in \cite{MR2829871}. To formulate it, Let $W_1, \cdots, W_N$ be $N$ independent, centered subgaussian variables with variance $\sigma$ and set $x_1,\ldots,x_N$ to be $N$  deterministic vectors in $\R^d$. Assume that ``design matrix", $\Gamma=N^{-1/2} \sum_{i=1}^N \inr{x_i,\cdot}e_i$, whose rows are $x_i/\sqrt{N}$, satisfies some Restricted Isometry Property (cf. \cite{MR2236170}). If $Y_i=\inr{x_i, t^*} + W_i, \ i=1,\ldots,N$, then for a well chosen regularization parameter $\lambda$, one has, with high probability,
\begin{equation}\label{eq:vlad_lasso}
\E\inr{X,\hat t-t^*}^2 \leq C \min\left\{\frac{\sigma^2\norm{t^*}_0 \log d}{n}, \sigma\norm{t^*}_1\sqrt{\frac{\log d}{n}}\right\}
\end{equation}
for a suitable absolute constant $C$. 

The error rate from \eqref{eq:vlad_lasso} consists of two components: the sparsity-dependent error term $\sigma^2(\norm{t^*}_0 \log d)/n$, and the complexity-dependent error term $\sigma \norm{t^*}_1\sqrt{(\log d)/n}$, and in what follows we will present a few other examples that combine the two rates -- because the procedure one uses to obtain both types of rate is the same.

The aim of this article is to address the ``complexity-based" aspect of the problem: to study regularization problems in which one believes that the $\Psi(f^*)$ is relatively small, and obtain an error rate that depends on $\Psi(f^*)$ rather than on some sparsity property of $f^*$.

\subsection{Attaining Minimax rates} \label{sub:minimax_rates_and_aims_of_the_paper}
A natural benchmark for measuring the success of a regularization method is the minimax error rate, assuming that $\Psi(f^*)$ is known: if one is given additional information on $\Psi(f^*)$, e.g., that $f^* \in \{ f : \Psi(f) \leq R\}$, one may consider the estimation problem in $\{ f : \Psi(f) \leq R\}$ using the given random data. Such a problem has an optimal error rate (called the minimax rate): it is the best rate any learning procedure may achieve in the class $\{ f : \Psi(f) \leq R\}$ given the random data $(X_i,Y_i)_{i=1}^N$. This minimax rate will serve as our benchmark, and will be compared with the error rates that we obtain.

Of course, one {\it is not} given additional information on $\Psi(f^*)$ and it is reasonable to expect that the error rate of the regularization procedure will be significantly slower than this  benchmark. The question we shall study here focuses on that gap. In fact, we will show that the price one has to pay for not knowing $\Psi(f^*)$ is surprisingly small, under rather weak assumptions.

\vskip0.4cm

From a technical perspective, all regularization-based procedures share one crucial aspect: the calibration of
the regularization parameter $\lambda$. That choice is very important as $\lambda$ is an essential component in ensuring that the error rate of the estimator $\hat{f}$ is well-behaved. Thus, to study the gap between the regularization error rate and the minimax rate, one has to identify the right choice of $\lambda$.

\begin{Question} \label{qu:reg}
What is the `correct choice' of the regularization parameter $\lambda$, and given that choice, what is the rate of convergence of RERM? Specifically, how far is the resulting rate from the one that could have been achieved had $\Psi(f^*)$ been given in advance?
\end{Question}
An answer to Question \ref{qu:reg} requires one to identify $\lambda$; to find a high probability upper bound on $\|\hat{f}-f^*\|_{L_2(\mu)}^2$ for that choice of $\lambda$; and then to compare the error rate to the minimax rate of the estimation problem in the ``true model'' $\{f : \Psi(f) \leq \Psi(f^*)\}$.

The strategy we use below follows a similar path to \cite{LM_sparsity} and is based on the small ball method, introduced in \cite{Shahar-ACM, Shahar-Gelfand,Shahar-Vladimir,shahar_general_loss}.

\subsection{The small-ball method} \label{sec:small-ball}
Given a closed and convex class $F$ and an unknown target $Y$, recall that $f^* \in F$ is a minimizer in $F$ of the functional $f \to \E(f(X)-Y)^2$.

The excess loss functional associated with $f \in L_2(\mu)$ is
\begin{align} \label{eq:basic-loss}
f \to {\cal L}_f(X,Y) = & \ell_f(X,Y)-\ell_{f^*}(X,Y) = (f(X)-Y)^2-(f^*(X)-Y)^2 \nonumber
\\
= & (f-f^*)^2(X)+2(f^*(X)-Y)(f-f^*)(X).
\end{align}
Moreover, since $F$ is closed and convex, then by the characterization of the nearest point map in a Hilbert space,
$$
\E (f^*(X)-Y)(f-f^*)(X) \geq 0 \ \ {\rm for \ every \ } f \in F;
$$
thus
\begin{align} \label{eq:basic-loss-1}
 \frac{1}{N}\sum_{i=1}^N (f^*(X_i)-Y_i)(f-f^*)(X_i)
\geq  \frac{1}{N}\sum_{i=1}^N (f^*(X_i)-Y_i)(f-f^*)(X_i) - \E (f^*(X)-Y)(f-f^*)(X).
\end{align}

Let $E$ be a subspace that contains $F$ and set $\Psi(\cdot)$ to be a regularization function on $E$ (i.e., a functional that satisfies Assumption~\ref{assum:reg-function}). Set $\rho \geq 0$ and put
$$
K_\rho(f^*) =\{h \in E : \Psi(h-f^*) \leq \rho\},
$$
which, by the convexity of $\Psi$, is a convex set.

\begin{Definition} \label{def:lambda-excess-loss}
For every $\lambda>0$ and any $f\in L_2(\mu)$, define the \textit{regularized excess loss} by
$$
{\cal L}_f^\lambda = \left(\ell_f+\lambda \Psi(f)\right)-\left(
\ell_{f^*}+\lambda \Psi(f^*)\right) =  {\cal L}_f + \lambda \left(\Psi(f)-\Psi(f^*)\right).
$$
\end{Definition}

Note that for every sample $(X_i,Y_i)_{i=1}^N$, a minimizer  $\hat{f}$ of the empirical regularized loss functional \eqref{eq:RERM} also minimizes in $F$ the empirical regularized excess loss $f\to P_N {\cal L}_f^\lambda$. Hence, since ${\cal L}_{f^*}^\lambda=0$, it follows that for every $(X_i,Y_i)_{i=1}^N$, the empirical regularized excess loss in $\hat f$ is non-positive:
\begin{equation}
  \label{eq:1}
 P_N {\cal L}_{\hat{f}}^\lambda \leq 0.
\end{equation}
This observation is at the heart of our analysis, as it allows one to
exclude functions $f$ in $F$ that satisfy $P_N {\cal L}_{f}^\lambda
> 0$ as potential minimizers of the empirical regularized loss function. Our strategy is therefore to show that if $f\in F$ and $\norm{f-f^*}_{L_2(\mu)}$ is not `too small', then necessarily $P_N {\cal L}_{f}^\lambda >0$ (for the right choice of $\lambda$); hence, functions cannot be minimizers of the empirical regularized (excess) loss function.

\vskip0.4cm
To simplify notation, set $\xi = Y - f^*(X)$,
$$
{\cal M}_{f -f^*}(X,Y)=\xi(f-f^*)(X)- \E \xi(f-f^*)(X)  \ \ {\rm and} \ \ {\cal Q}_{f-f^*}(X)=(f-f^*)^2(X);
$$
therefore, combining \eqref{eq:basic-loss} and \eqref{eq:basic-loss-1},
\begin{equation}
  \label{eq:decomp}
P_N {\cal L}_f \geq P_N {\cal Q}_{f-f^*} -2\left|P_N {\cal M}_{f-f^*}\right|.
\end{equation}

The main step in the small-ball method is to find a lower bound on the quadratic process $f \to P_N {\cal Q}_{f-f^*}$ and an upper bound on $f \to \left|P_N {\cal M}_{f-f^*}\right|$. The two estimates should hold with high probability on certain subsets of $F$. Then, they have to be compared with the behaviour of the regularization term $\lambda(\Psi(f)-\Psi(f^*))$ on those sets to ensure that $P_N{\cal L}_f^\lambda>0$.

\vskip0.4cm

A uniform lower bound on the quadratic component $P_N {\cal Q}_{f-f^*}$ can be obtained under a weak assumption called the {\it small-ball} condition:
\begin{Assumption} \label{ass:small-ball}
Assume that there are constants $\kappa>0$ and $0<\eps \leq 1$, for which, for every $f,h \in F$,
$$
Pr\left(|f-h| \geq \kappa \|f-h\|_{L_2(\mu)}\right) \geq \eps.
$$
\end{Assumption}
There are numerous examples in which Assumption \ref{ass:small-ball} may be
verified for $\kappa$ and $\eps$ that are absolute constants and we refer the reader to
\cite{Shahar-Gelfand,Shahar-ACM,LM_compressed,shahar_general_loss,Shahar-Vladimir,RV_small_ball}
for some of them.

To put assumption~\ref{ass:small-ball} on $X$ in some perspective, recall that the class $F=\{f_t=\inr{\cdot,t}:t\in\R^d\}$ is {\it identifiable} if for every
$t_1,t_2\in\R^d$, $Pr(f_{t_1} \neq f_{t_2})>0$, (where the probability is taken with respect to the underlying measure $\mu$). By linearity, this condition is equivalent to assuming that for every $t \in \R^d$, $Pr(|\inr{X,t}|>0)>0$. Thus, the small-ball condition is simply a uniform estimate on the degree of identifiability of class $F$ and is therefore a rather weak assumption.

\vskip0.5cm
Now, let us introduce two complexity parameters that play a central role in our analysis. Let $D$ be the unit ball in $L_2(\mu)$ and for $r>0$ set
$$
rD_{f^*}=\{f \in L_2(\mu): \|f-f^*\|_{L_2(\mu)} \leq r\}=f^*+rD.
$$

\begin{Definition} \label{def:r-Q}
Given a class $F$ of functions and $\tau>0$, let
$$
r_Q(F,\tau) = r_Q(F, f^*, \tau) = \inf\left\{r>0: \E \sup_{f \in F \cap rD_{f^*}} \left|\frac{1}{N}\sum_{i=1}^N \eps_i (f-f^*)(X_i)\right| \leq \tau r \right\},
$$
where $(\eps_i)_{i=1}^N$ are independent, symmetric, $\{-1,1\}$-valued random variables that are also independent of $(X_i,Y_i)_{i=1}^N$.

Set
\begin{equation} \label{eq:phi}
\phi_N(F,f^*,s)= \sup_{f \in F \cap sD_{f^*}} \left|\frac{1}{\sqrt{N}} \sum_{i=1}^N \eps_i \xi_i (f-f^*)(X_i)\right|
\end{equation}
and put
$$
r_M(F,\tau,\delta) = r_M(F,f^*,\tau,\delta) =  \inf\left\{ s>0 : Pr \left(  \phi_N(F,f^*,s)\leq \tau s^2 \sqrt{N} \right) \geq 1-\delta \right\}.
$$
\end{Definition}

One may show the following (see Theorem~3.1 in \cite{Shahar-ACM}):
\begin{Theorem} \label{thm:small-ball-method}
Let $F$ be a closed, convex class of functions that satisfies
Assumption \ref{ass:small-ball} with constants $\kappa$ and $\eps$, and
set $\theta = \kappa^2 \eps/16$. For every $\delta\in(0,1)$, with probability at least $1-\delta-2\exp(-N\eps^2/2)$ one has both:
\begin{description}
\item{$\bullet$} for every $f \in F$,
$$
|P_N {\cal M}_{f-f^*}| \leq \frac{\theta}{4}\max\left\{\|f-f^*\|_{L_2(\mu)}^2,r_M^2 \left(F,\theta/5,\delta/4\right)\right\},
$$
\item{$\bullet$} for every $f \in F$ with $\|f-f^*\|_{L_2(\mu)} \geq r_Q \left(F,\kappa \eps/32\right)$,
$$
P_N {\cal Q}_{f-f^*} \geq \theta \|f-f^*\|_{L_2(\mu)}^2.
$$
\end{description}
In particular, with probability at least $1-\delta-2\exp(-N\eps^2/2)$, $P_N {\cal L}_f \geq \frac{\theta}{2}
  \|f-f^*\|_{L_2(\mu)}^2$ for every $f\in F$ that satisfies
$$
\|f-f^*\|_{L_2(\mu)} \geq \max\left\{r_M \left(F,\theta/5,\delta/4\right),r_Q \left(F,\kappa \eps/32\right)\right\}.
$$
\end{Theorem}
\vskip0.4cm

\begin{Remark}
An immediate outcome of Theorem \ref{thm:small-ball-method} is that with high probability, a minimizer in $F$ of the empirical excess-loss functional $P_N {\cal L}_f$ must satisfy
\begin{equation} \label{eq-ERM-est}
\|\tilde{f}-f^*\|_{L_2(\mu)} \leq \max\left\{r_M \left(F,\theta/5,\delta/4\right),r_Q \left(F,\kappa \eps/32\right)\right\}.
\end{equation}
\end{Remark}

In fact, results from \cite{LM13} show that \eqref{eq-ERM-est} is optimal in the minimax sense under additional mild technical assumptions on $F$ when the data are assumed to satisfied the Gaussian regression model, that is, when the targets are of the form $Y=f_0(X)+W$ for $f_0 \in F$ and $W$ that is a centered Gaussian random variable, independent of $X$. {\it Empirical risk minimization} performed in the set
$$
F^*=\{f \in F : \Psi(f) \leq \Psi(f^*)\}
$$
yields
\begin{equation} \label{eq:error-estimate}
\|\tilde{f}-f^*\|_{L_2} \leq \max\left\{r_M \left(F^*,\theta/5,\delta/4\right),r_Q \left(F^*,\kappa \eps/32\right)\right\},
\end{equation}
and the r.h.s. of \eqref{eq:error-estimate} is the minimax rate of the estimation problem in $F^*$ (up to the technical assumptions mentioned earlier); it will serve as a benchmark for the performance of the regularization procedure \eqref{eq:RERM}.

\subsection{The Main result} \label{sec:main-result}
Let $F \cap K_\rho(f^*)=\{f \in F : \Psi(f-f^*) \leq \rho\}$ and observe that these are convex subsets of $F$. To simplify notation, set
\begin{equation} \label{eq:r-M-r-Q-choice}
r_M(\rho)=r_M \Big(F \cap K_{\rho}(f^*), \frac{\kappa^2\eps}{80},\frac{\delta}{4}\Big) \ \ {\rm and} \ \ r_Q(\rho)=r_Q \Big(F \cap K_\rho(f^*),\frac{\kappa \eps} {32}\Big),
\end{equation}
and let $r(\cdot)$ be a function that satisfies for every $\rho \geq 0$
\begin{equation}
  \label{eq:r-rho}
 r(\rho)\geq\max\{r_M(\rho),r_Q(\rho)\}.
\end{equation}
It should be noted that $r(\rho)$ may depend on $f^*$, and that it does depend on other parameters -- like $\delta$, $\kappa$ and $\eps$. We will not specify the dependence on those parameters, but rather, only on the radius $\rho$.

The geometry of the sets $F \cap K_\rho(f^*)$ (see Figure~\ref{fig:partition_set_F}) determine both the error rate and the regularization parameter $\lambda$, and $r(\rho)$ measures the sets' `sizes'.

The choice of $\lambda$ is made as follows:

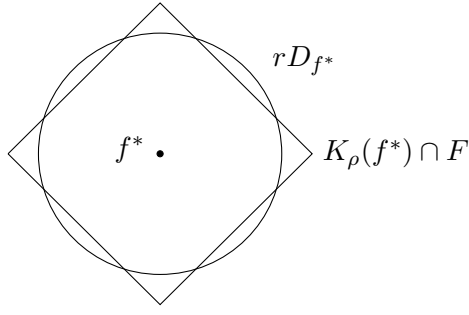
\begin{figure}[!h]
\centering
\begin{tikzpicture}[scale=0.2]
\draw (0,0) circle (8cm);
\draw (-2,0.3) node {$f^*$};
\filldraw (0,0) circle (0.2cm);
\draw (-10,0) -- (0,10) -- (10,0) -- (0,-10) -- (-10,0);
\draw (15.5,0) node {$K_{\rho}(f^*)\cap F$};
\draw (9.5,6.2) node {$r D_{f^*}$};
\end{tikzpicture}
\caption{Localization of the set $F\cap K_\rho(f^*)$, i.e. its intersection with $L_2(\mu)$-balls of various radii $r$ for the right choice of radius $\rho$, plays a central role in the analysis of the quadratic and multiplier processes.}
\label{fig:partition_set_F}
\end{figure}

Let
$$
{\cal O}(\rho) = \sup \Big(|P_N {\cal M}_{f-f^*}| : f \in F \cap K_{\rho}(f^*) \cap r(\rho)D_{f^*}\Big)
$$
and for $\tau>0$ and $0<\delta<1$, set
$$
\gamma_{\cal O}(\rho, \tau, \delta) = \inf\left\{x>0: Pr\left({\cal O}(\rho) \leq \tau x \right) \geq 1-\delta \right\}
$$
and 
$$
\gamma_{\cal O}(\rho) = \gamma_{\cal O}(\rho, 3/(80\eta^3), \delta).
$$
In other words, $\gamma_{\cal O}(\rho)$ is proportional to the smallest possible upper estimate on ${\cal O}(\rho)$ that still holds with probability at least $1-\delta$.

\begin{Definition}  \label{def:lambda0}
For any  $\tau>0$ and $0<\delta<1$, set
\begin{equation*}
  \lambda_0(\delta,\tau)=\sup_{\rho>0, f^* \in F}\frac{\gamma_{\cal O}(\rho, \tau, \delta)}{\rho}.
\end{equation*}
\end{Definition}

To compare $\lambda_0(\delta,\tau)$  with $r_Q$ and $r_M$, first note that $r_M(\rho)$ and $\cO(\rho)$ both depend on properties of the multiplier processes indexed by localizations of $F\cap K_\rho(f^*)$, and recall that symmetrized and centered processes are essentially equivalent. Second, if $r(\rho)=r_M(\rho)$ then $\gamma_{\cal O}(\rho) \sim r_M^2(\rho)$; moreover, $\gamma_{\cal O}(\rho)$ is trivially bounded by $\sim r^2(\rho)$ for the right choice of $\tau$ and $\delta$. However, if $r_M(\rho) \leq r_Q(\rho)$, that is, when $r(\rho)=r_Q(\rho)$   -- which is the case when $\rho$ is very large -- one may find that $\gamma_{\cal O}(\rho)$ is actually significantly smaller than $\sim r^2(\rho)$. This observation is of crucial importance because of the choice of the regularization parameter: for the right choice of $\tau$, $\gamma_{{\cal O}}(\rho) \leq r^2(\rho)$ and
$$
\lambda_0(\delta,\tau) \leq \sup_{\rho>0, f^* \in F} r^2(\rho)/\rho;
$$
thus, one may be tempted to select the latter as a regularization term. However, there are natural examples in which $\sup_{\rho>0} r^2(\rho)/\rho = \infty$, rendering that choice impossible, whereas $\sup_{\rho>0, f^* \in F}{\gamma_{\cal O}(\rho)}/{\rho}$ turns out to be finite. Of course, there are still cases in which $\sup_{\rho>0, f^* \in F} r^2(\rho)/\rho$ is finite, and $\lambda_0(\delta,\tau)$ is of the same order as $\sup_{\rho>0, f^* \in F} r_M^2(\rho)/\rho$, though that is not the generic situation.

\vskip0.4cm

We now come to the main result of the article.
\begin{Theorem}\label{thm:main}
Let $F$ be a closed, convex class of functions that satisfies Assumption \ref{ass:small-ball} with constants $\kappa$ and $\eps$. Set $\Psi(\cdot)$ to be a regularization function that satisfies Assumption~\ref{assum:reg-function} with constant $\eta$. Furthermore, assume that $\lim_{\rho  \to 0} r(\rho)= 0$ and put $\lambda> \lambda_0(\delta, 3/(80\eta^3))$.

If $\hat{f}$ is the RERM with a regularization parameter $\lambda$ as in \eqref{eq:RERM}, then with probability at least $1-2\delta-2\exp(-N\eps^2/2)$,
\begin{equation} \label{eq:in-main-form}
\|\hat{f}-f^*\|_{L_2(\mu)}^2 \leq \max \Bigl\{r^2\big(10\eta \Psi(f^*)\big),  \Big(\frac{32}{\kappa^2\eps}\Big)\lambda \Psi(f^*)\Bigr\}.
\end{equation}
\end{Theorem}

Observe that $\lambda_0$ depends only on the oscillations of the multiplier process. Hence, if the problem is noise-free then $\lambda_0=0$, showing that any regularization parameter $\lambda>0$ would do. Moreover, in that case $r_M(\rho)=0$ and so one can choose $r(\rho)\geq r_Q(\rho)$ obtaining an error rate that depends only on $r_Q^2(10 \eta \Psi(f^*))$.

\vskip0.4cm

As noted previously, if one considers empirical risk minimization performed in $F^*=\{f \in F : \Psi(f) \leq \Psi(f^*)\}$, the resulting error rate is $\|\tilde{f}-f^*\|_{L_2(\mu)}^2 \leq c_0 r^2(c\Psi(f^*))$ for a suitable absolute constant $c_0$ and a constant $c$ that depends on $\kappa$, $\eps$ and $\delta$; moreover, under some minor additional assumptions, that rate is optimal in the minimax sense (cf. \cite{LM13}) when one takes $r(\rho) \sim \max\big\{r_M(\rho), r_Q(\rho)\big\}$. Hence, up to constants involved, the first term in Theorem \ref{thm:main} is essentially the minimax rate that one can obtain if $\Psi(f^*)$ were known.

If one chooses $\lambda\sim \lambda_0(\delta, \tau)$ for $\tau = 3/(80 \eta^3)$ then the second term in \eqref{eq:in-main-form} is of the order of
$$
\lambda \Psi(f^*) = \Big(\sup_{\rho,f^*} \frac{\gamma_{\cal O}(\rho, \tau, \delta)}{\rho}\Big) \cdot \Psi(f^*).
$$
Note that for $\rho$ that is of the order of $\Psi(f^*)$, one has
$$
\frac{\gamma_{\cal O}(\rho, \tau, \delta)}{\rho} \cdot \Psi(f^*) \leq c_1\gamma_{\cal O}(\rho, \tau, \delta) \leq c_2r^2(c_3\Psi(f^*)),
$$
which coincides with the first term, up to the constants involved. Thus, the price that one has to pay for not knowing $\Psi(f^*)$ is manifested in the need to take the supremum over all possible choices of $\rho$ in the second term, rather than considering only the level $\rho \sim \Psi(f^*)$ of the ``true model".

Thankfully, there are many natural cases in which that price is rather small, allowing for satisfactory outcomes of Theorem \ref{thm:main} closed to the minimax rate.

\vskip0.4cm

We end this introduction with a word about notation. Throughout, absolute constants or constants that depend on other parameters are denoted by $c$, $C$, $c_1$,
$c_2$, etc., (and, of course, we will specify when a constant is
absolute and when it depends on other parameters). The values of these
constants may change from line to line. The notation $x\sim y$ (resp. $x\lesssim y$) means that there exist absolute constants $0<c<C$ for which $cy\leq x\leq Cy$ (resp. $x\leq Cy$). If $b>0$ is a parameter, then $x\lesssim_b y$ means that $x\leq C(b) y $ for some constant $C(b)$ that depends only on $b$.

The normed space $\ell_p^d$ is $\R^d$ endowed with the norm $\|x\|_{p}=\big(\sum_{j}|x_j|^p\big)^{1/p}$; the corresponding unit ball is denoted by $B_p^d$ and the unit Euclidean sphere in $\R^d$ is $S^{d-1}$.

Finally, from here on we will write $Pr$ and $\| \ \|_{L_2}$ without specifying the underlying measure.

\section{Proof of Theorem \ref{thm:main}} \label{sec:proof}
The proof of Theorem \ref{thm:main} follows an almost identical path as the proof of Theorem 3.2 from \cite{LM_sparsity}. The differences in the two arguments are minor and their source is the fact that unlike \cite{LM_sparsity}, here we do not assume that $\Psi$ is a norm. We will outline in Remark \ref{rem:proof-of-main} how a version of Theorem \ref{thm:main} may be derived directly from Theorem 3.2 in \cite{LM_sparsity} when $\Psi$ is a norm.

\vskip0.4cm

Theorem \ref{thm:main} is an immediate outcome of the following lemma:
\begin{Lemma} \label{lemma:main}
Let $\lambda_0=\lambda_0(\delta,3/(80\eta^3))$ and set
$\lambda > \lambda_0$. If $\lim_{\rho \to 0} r(\rho)=
0$, $\rho \geq 5\eta \Psi(f^*)$ and $\rho>0$, then with probability at least $1-2\delta-2\exp(-N\eps^2/2)$,
$$
\|\hat{f}-f^*\|_{L_2}^2 \leq \max\left\{r^2(\rho),(32/(\kappa^2\eps))\lambda\Psi(f^*)\right\}.
$$
\end{Lemma}

To see how Lemma \ref{lemma:main} can be used to conclude the proof of Theorem
\ref{thm:main}, observe that if $\Psi(f^*)>0$, one may simply select
$\rho=5\eta\Psi(f^*)$ in the lemma.  If, on the other hand,
$\Psi(f^*)=0$, let $(\gamma_n)_{n=1}^\infty$ be a positive sequence
decreasing to $0$ and set ${\cal A}_n=\{\|\hat{f}-f^*\|_{L_2} \leq
\gamma_n\}$, which is a decreasing sequence of events. If
$Pr({\cal A}_n) \geq 1-\nu$ for some $0 \leq \nu\leq 1$ and every $n$ then $Pr(\{\hat{f}=f^*\})
\geq 1-\nu$. Since $\lim_{\rho \to 0} r(\rho)= 0$, one may apply Lemma
\ref{lemma:main} to each member of a nonnegative sequence $\rho_n$ that decreases to zero and for which $\gamma_n =r(\rho_n)$ decreases to zero. By Lemma~\ref{lemma:main}, $Pr({\cal A}_n) \geq 1-2\delta-2\exp(-N\eps^2/2)$ for every $n$ and the proof of Theorem \ref{thm:main} follows.

\vskip0.5cm

\noindent{\bf Proof of Lemma \ref{lemma:main}.} Fix $f^*$ and set $\rho>0$ that satisfies $\rho \geq 5\eta \Psi(f^*)$.
Let
$$
F_1=\{f \in F : \Psi(f-f^*) \leq \rho\} = F \cap K_{\rho}(f^*),
$$
and
$$
F_2=\{f \in F : \Psi(f-f^*) = \rho\}.
$$
Clearly, $F_1$ is a convex set that contains $f^*$, and by the continuity of the real-valued function $t \to \Psi(f^*+t (f-f^*))$, every ray $[f^*,f)$ that originates in $f^*$ and passes through some $f \in F \backslash F_1$ intersects $F_2$.

Let $\theta=\kappa^2\eps/16$ and set
$$
    r_Q(\rho)=r_Q(F_1,\kappa\eps/32) \ \ {\rm and} \ \  r_M(\rho)=r_M(F_1,\theta/5,\delta/4).
$$
There is an event ${\cal A}_0$ of probability at least $1-\delta-2\exp(-N\eps^2/2)$,
and for every $(X_i,Y_i)_{i=1}^N \in {\cal A}_0$ the following holds:
\begin{description}
\item{$\bullet$} If $f \in F_1$ and $\|f-f^*\|_{L_2} \geq r_Q(\rho)$ then
$$
\frac{1}{N}\sum_{i=1}^N (f-f^*)^2(X_i) \geq \theta \|f-f^*\|_{L_2}^2.
$$
\item{$\bullet$} If $f \in F_1$ then
$$
\left|\frac{1}{N}\sum_{i=1}^N \xi_i(f-f^*)(X_i) - \E \xi (f-f^*)(X)\right| \leq \frac{\theta}{4} \max\{\|f-f^*\|_{L_2}^2,r_M^2(\rho)\}.
$$
\end{description}

In particular, if $f \in F_1$ and $\|f-f^*\|_{L_2}
  \geq r(\rho)\geq\max\{r_M(\rho), r_Q(\rho)\}$ then
$$
P_N {\cal L}_f \geq \frac{\theta}{2} \|f-f^*\|_{L_2}^2.
$$
By the choice of $\lambda$, there is an event ${\cal A}_1$ of probability at least $1-\delta$ on which if $f\in F_1$ and
$\norm{f-f^*}_{L_2}\leq r(\rho)$, then
\begin{equation}
\label{eq:multi_lambda0}
 \left| \frac{2}{N}\sum_{i=1}^N \xi_i(f-f^*)(X_i) - \E \xi (f-f^*)(X) \right| < \frac{3}{80\eta^3}\lambda
  \rho<\frac{3}{5\eta}\lambda \rho.
\end{equation}

Set ${\cal A}= {\cal A}_0\cap{\cal A}_1$ and let $(X_i,Y_i)_{i=1}^N \in {\cal A}$. The proof now follows in three steps:
\begin{description}
\item{(1)} Show that the functional $f \to P_N {\cal L}_f^\lambda$ is bounded from below -- away from zero -- in $F_2$.
\item{(2)} An outcome of ($1$) is that if $f \in F
  \backslash F_1$, $P_N {\cal L}_f^\lambda >0$; hence, $\hat{f} \not \in F \backslash F_1$.
\item{(3)} Finally, pin-point $\hat{f}$ within $F_1 = \{f \in F : \Psi(f-f^*) \leq \rho\}$.
\end{description}

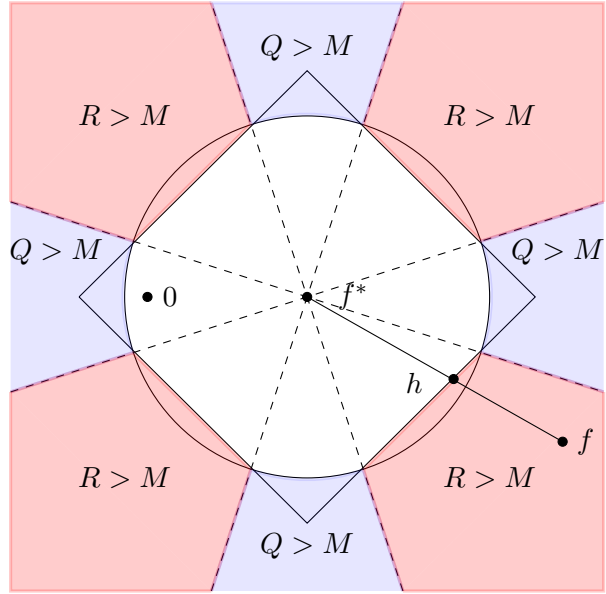
\begin{figure}[!h]
\centering
\begin{tikzpicture}[scale=0.3]
\draw (0,0) circle (8cm);
\draw (2,0) node {$f^*$};
\filldraw (0,0) circle (0.2cm);
\draw (-6,0) node {$0$};
\filldraw (-7,0) circle (0.2cm);
\draw (-10,0) -- (0,10) -- (10,0) -- (0,-10) -- (-10,0);
\draw[style = dashed] (-13,-4.2) -- (13,4.2);
\draw[style = dashed] (-13,4.2) -- (13,-4.2);
\draw[style = dashed] (-4.2,-13) -- (4.2,13);
\draw[style = dashed] (4.2,-13) -- (-4.2,13);
\filldraw[color=red, very thick, fill = red, opacity = 0.2] (4.2,13) -- (2.47,7.64) -- (7.64,2.47) -- (13,4.2);
\filldraw[color=red, very thick, fill = red, opacity = 0.2] (4.2,13) -- (13,13) -- (13,4.2);
\filldraw[color=red, very thick, fill = red, opacity = 0.2] (-4.2,-13) -- (-2.47,-7.64) -- (-7.64,-2.47) -- (-13,-4.2);
\filldraw[color=red, very thick, fill = red, opacity = 0.2] (-4.2,-13) -- (-13,-13) -- (-13,-4.2);
\filldraw[color=red, very thick, fill = red, opacity = 0.2] (-4.2,13) -- (-2.47,7.64) -- (-7.64,2.47) -- (-13,4.2);
\filldraw[color=red, very thick, fill = red, opacity = 0.2] (-4.2,13) -- (-13,13) -- (-13,4.2);
\filldraw[color=red, very thick, fill = red, opacity = 0.2] (4.2,-13) -- (2.47,-7.64) -- (7.64,-2.47) -- (13,-4.2);
\filldraw[color=red, very thick, fill = red, opacity = 0.2] (4.2,-13) -- (13,-13) -- (13,-4.2);
\filldraw[color=blue, very thick, fill = blue, opacity = 0.1] (-4.2,13) -- (-2.47,7.64) arc (110:70:7.2) -- (4.2,13) -- cycle;
\filldraw[color=blue, very thick, fill = blue, opacity = 0.1] (-4.2,-13) -- (-2.47,-7.64) arc (-110:-70:7.2) -- (4.2,-13);
\filldraw[color=blue, very thick, fill = blue, opacity = 0.1] (-13,4.2) -- (-7.64,2.47) arc (160:200:7.2) -- (-13,-4.2);
\filldraw[color=blue, very thick, fill = blue, opacity = 0.1] (13,4.2) -- (7.64,2.47) arc (20:-20:7.2) -- (13,-4.2);
\draw (8,8) node {$R>M$};
\draw (-8,8) node {$R>M$};
\draw (8,-8) node {$R>M$};
\draw (-8,-8) node {$R>M$};
\draw (0,11) node {$Q>M$};
\draw (-11,2) node {$Q>M$};
\draw (11,2) node {$Q>M$};
\draw (0,-11) node {$Q>M$};
\draw (0,0) -- (11.2,-6.4);
\draw (12.2,-6.4) node {$f$};
\filldraw (11.2,-6.4) circle (0.2cm);
\draw (4.7,-3.8) node {$h$};
\filldraw (6.42,-3.63) circle (0.2cm);
\end{tikzpicture}
\caption{$P_N {\cal L}_f^\lambda>0$ for two different reasons:  either $Q>M$ -- the quadratic component dominates the multiplier component, or $R>M$ -- the regularization component dominates the multiplier component. Unlike Theorem~3.2 in \cite{LM_sparsity}, here we choose $\rho\sim \Psi(f^*)$ to ensure that $0\in F \cap K_{\rho}(f^*)$.}
\label{fig:partition_set_F_positive_excess_loss}
\end{figure}

\noindent{\bf Step 1.} Fix $f \in F_2$ and note that by the `triangle inequality' satisfied by $\Psi$,
$$
\Psi(f) \geq \eta^{-1}\Psi(f-f^*)-\Psi(f^*).
$$
Recall that $\eta^{-1}\Psi(f-f^*) \geq \eta^{-1}\rho \geq 5 \Psi(f^*)$ and thus,  $\Psi(f)-\Psi(f^*) \geq (3/5)\eta^{-1}\rho$. Hence, if $\|f-f^*\|_{L_2} \geq r(\rho)$ then
$$
P_N {\cal L}_f^\lambda \geq (\theta/2) \|f-f^*\|_{L_2}^2 + \lambda \rho \cdot \frac{3}{5\eta} > 0.
$$
On the other hand, by the choice of $\lambda$, if $\|f-f^*\|_{L_2} \leq r(\rho)$ then
\begin{align*}
P_N {\cal L}_f^\lambda \geq & -\left|\frac{1}{N}\sum_{i=1}^N \xi_i (f-f^*)(X_i) - \E \xi (f-f^*)(X)\right|+ \lambda \left(\Psi(f)-\Psi(f^*)\right)
\\
\geq & -\left|\frac{1}{N}\sum_{i=1}^N \xi_i (f-f^*)(X_i) - \E \xi (f-f^*)(X)\right| + \lambda \rho \cdot \frac{3}{5\eta} >0.
\end{align*}

It should be noted that the same proof shows that on the event ${\cal A}$, for every $f \in F_2$,
\begin{equation} \label{eq:in-F_2}
P_N {\cal L}_f + \frac{\lambda}{2\eta^2}(\Psi(f)-\Psi(f^*))>0,
\end{equation}
a fact that will be used below. Indeed, $(\lambda/2\eta^2) \cdot (\Psi(f)-\Psi(f^*)) \geq (\lambda/2\eta^2) \cdot (3\rho/5\eta)$ and by \eqref{eq:multi_lambda0}, if $\|f-f^*\|_{L_2} \leq r(\rho)$ then $P_N{\cal L}_f \geq -(3/80) \cdot (\lambda \rho/\eta^3)$.
\vskip0.3cm

\noindent{\bf Step 2.} Let $f \in F \backslash F_1$ and note that by the convexity of $F$ and the continuity of $\Psi$ on rays, there is some $h \in F_2$ and $R > 1$ for which $f=f^*+R(h-f^*)$. Thus,
$$
P_N {\cal L}_f^\lambda = \frac{R^2}{N} \sum_{i=1}^N (h-f^*)^2(X_i) + \frac{2R}{N}\sum_{i=1}^N \xi_i(h-f^*)(X_i) + \lambda \left(\Psi(f)-\Psi(f^*)\right).
$$
Observe that
\begin{equation} \label{eq:difference-in-proof}
\Psi(f)-\Psi(f^*) \geq \frac{R}{2\eta^2} \left(\Psi(h)-\Psi(f^*)\right);
\end{equation}
indeed,
$$
\Psi(f^*+R(h-f^*)) \geq \eta^{-1}\Psi(R(h-f^*)) - \Psi(f^*) \geq R\eta^{-1}\Psi(h-f^*) - \Psi(f^*),
$$
and thus it suffices to show that
$$
\frac{R}{\eta}\Psi(h-f^*) \geq \frac{R}{2\eta^2}\Psi(h)+2\Psi(f^*).
$$
But since $\Psi(h-f^*) \geq 5\eta \Psi(f^*)$, $\eta \geq 1$ and $R \geq 1$, one has
\begin{align*}
&\frac{R}{\eta}\Psi(h-f^*) \geq  \frac{R}{2\eta}\Psi(h-f^*)+\frac{R}{2}\Psi(f^*)+2R\Psi(f^*)
\\
\geq & \frac{R}{2\eta}\left(\Psi(h-f^*)+\Psi(f^*)\right)+2\Psi(f^*)
\geq  \frac{R}{2\eta^2}\Psi(h)+2\Psi(f^*),
\end{align*}
and \eqref{eq:difference-in-proof} follows.

\vskip0.3cm

Finally, applying \eqref{eq:in-F_2} to $h\in F_2$,
\begin{align*}
P_N {\cal L}_f^\lambda \geq & \frac{R^2}{N} \sum_{i=1}^N (h-f^*)^2(X_i) + \frac{2R}{N}\sum_{i=1}^N \xi_i(h-f^*)(X_i) + \lambda \frac{R}{2\eta^2} \left(\Psi(h)-\Psi(f^*)\right)
\\
\geq & R\left(P_N {\cal L}_h + \frac{\lambda}{2\eta^2}(\Psi(h)-\Psi(f^*))\right)>0,
\end{align*}
and $\hat{f} \not \in F \backslash F_1$.
\vskip0.3cm
\noindent{\bf Step 3.}
Turning to $F_1 = \{f \in F : \Psi(f-f^*) \leq \rho\}=F \cap K_{\rho}(f^*)$, recall that if $f \in F_1$ and $\|f-f^*\|_{L_2} \geq r(\rho)$, then $P_N {\cal L}_f \geq (\theta/2) \|f-f^*\|_{L_2}^2$; hence, if $f$ is a potential minimizer and $\|f-f^*\|_{L_2} \geq r(\rho)$ then
\begin{align*}
0 \geq P_N {\cal L}_f^\lambda \geq & (\theta/2) \|f-f^*\|_{L_2}^2 + \lambda\left(\Psi(f)-\Psi(f^*)\right)
\geq  (\theta/2) \|f-f^*\|_{L_2}^2 - \lambda \Psi(f^*),
\end{align*}
and
$$
\|\hat{f}-f^*\|_{L_2}^2 \leq \frac{2 \lambda}{\theta} \Psi(f^*),
$$
as claimed.
\endproof

\begin{Remark} \label{rem:proof-of-main}
If $\Psi$ happens to be a norm (which is an assumption slightly stronger than Assumption~\ref{assum:reg-function}), one may apply Theorem 3.2 from \cite{LM_sparsity} directly. Indeed, and using the notation from \cite{LM_sparsity} if $\rho \gtrsim \Psi(f^*)$ then the set $K=\{f : \Psi(f-f^*) \leq \rho/20\}$ contains a $\Psi$-ball around $0$, and $\Gamma_{f^*}(\rho)$ -- the collection of norming functionals (i.e. the sub-gradient of $\Psi$) of any $h \in K$ -- is the entire unit ball in the dual space to $(E,\Psi)$. Recall that
$$
\Delta(\rho)=\inf_h \sup_{z^* \in \Gamma_{f^*}(\rho)} z^*(h-f^*),
$$
where the infimum is taken in the set
$$
\{h \in F : \ \Psi(h-f^*) = \rho \ {\rm and} \ \|h-f^*\|_{L_2} \leq r(\rho)\}.
$$
Since $\Gamma_{f^*}(\rho)$ is the entire dual unit ball, it follows that $\Delta(\rho)=\rho$, and Theorem 3.2 in \cite{LM_sparsity} may be applied. The desired version of Theorem~\ref{thm:main} now follows from Remark 3.3 in \cite{LM_sparsity}.
\end{Remark}

\section{Towards the examples - preliminary estimates} \label{sec:pre:estimates}
It is rather obvious that any implementation of Theorem \ref{thm:main} requires specific estimates on $r_M$, $r_Q$ and $\lambda_0$. This section is devoted to some preliminary facts that will play an instrumental part in establishing such estimates.

\vskip0.5cm

Our main interest is the study of upper bounds on the three processes used to define the parameter $r_M$, $r_Q$ and $\gamma_\cO$, and which have the following forms:
$$
(*)=\sup_{f \in F} \left|\sum_{i=1}^N \eps_i \xi_i f(X_i) \right| , \ \ (**)=\sup_{f \in F} \left|\sum_{i=1}^N \left(\xi_i f(X_i) - \E \xi f(X)\right) \right|\ \\ {\mbox{ and }}  \\ \E \sup_{f\in F}\Big|\sum_{i=1}^N \eps_i f(X_i)\Big|,
$$
where $X_1,...,X_N$ are independent and distributed according to the underlying measure $\mu$, $\xi_1,...,\xi_N$ are independent copies of $\xi \in L_q$ for some $q>2$ (though $(\xi_i)_{i=1}^N$ need not be independent of $(X_i)_{i=1}^N$), and $(\eps_i)_{i=1}^N$ are independent, symmetric $\{-1,1\}$-valued random variables that are independent of $(X_i)_{i=1}^N$ and $(\xi_i)_{i=1}^N$.

Standard symmetrization methods (see, e.g., \cite{MR757767,LT:91,vanderVaartWellner}) show that $(*)$ and $(**)$ are equivalent in expectation and in deviation (up to a slight restriction on the deviation parameter). We will present one example in which this symmetrization argument is carried out in full (Theorem \ref{theo:lasso-small-ball-weak-moments}), but in the other examples we will only consider the symmetrized case.

\subsection{Estimates for subgaussian classes} \label{sub:estimates_for_subgaussian_classes}
The first result is from  \cite{shahar_multi_pro}, under the assumption that $F$ is an $L$-subgaussian class of functions.

\begin{Definition} \label{def:subgaussian-class}
A class of functions $F \subset L_2(\mu)$ is $L$-subgaussian if for every $f,h \in F \cup \{0\}$ and every $u \geq 1$,
$$
Pr(|f(X)-h(X)| \geq Lu\|f-h\|_{L_2(\mu)}) \leq 2\exp(-u^2/2)
$$where $X$ is distributed according to $\mu$.
\end{Definition}

Let $F \subset L_2(\mu)$ and set $\{G_f : f \in F\}$ to be the centered, canonical Gaussian process indexed by
$F$ (i.e., the covariance operator of the process is $\E G_f G_g = \E f(X)g(X)$ for every $f,g\in F$). Put
\begin{equation}\label{eq:gaussian_mean_width}
\ell^*(F) = \E \sup_{f \in F} G_{f}, \ \ {\rm and } \ \ d_2(F) = \sup_{f \in F} \|f\|_{L_2(\mu)}.
\end{equation}

\begin{Theorem}[Corollary~1.10 in \cite{shahar_multi_pro}]
  \label{theo:multiplier_process_shahar}
  Let $X$ be distributed according to $\mu$, set $\xi\in L_q$ for some $q>2$ and assume that $F \subset L_2(\mu)$ is an $L$-subgaussian class. There are constants $c,c_0,c_1,c_2$ and $c_3$ that depend only on $q$, for which, for any $w,u>c$, with probability at least
  \begin{equation*}
    1-\frac{c_0\log^q N}{w^{q}N^{q/2-1}}-2\exp\left(-c_1 u^2(\ell^*(F)/d_2(F))^2\right),
  \end{equation*}

\begin{equation*}
    \sup_{f\in F}\left|\frac{1}{N}\sum_{i=1}^N \eps_i\xi_i f(X_i)\right|\leq c_2 Lw u\norm{\xi}_{L_q}\frac{\ell^*(F)}{\sqrt{N}}
  \end{equation*}
  and
  \begin{equation*}
    \sup_{f\in F}\left|\frac{1}{N}\sum_{i=1}^N \xi_i f(X_i)-\E \xi f(X)\right|\leq c_2 Lw u\norm{\xi}_{L_q}\frac{\ell^*(F)}{\sqrt{N}}.
  \end{equation*}
\end{Theorem}

\begin{Corollary}
  \label{coro:multiplier_process}
Using the notation and assumptions of Theorem \ref{theo:multiplier_process_shahar}, let $\xi = Y-f^*(X)$ and assume that $\xi \in L_q$ for some
$q>2$. Fix $\tau>0$ and $0<\delta<1$, and set $A>0$ for which
\begin{equation}
\label{eq:def_rM}
  c_2 L w u \norm{\xi}_{L_q} \ell^*\left(F\cap A D_{f^*}\right)\leq \tau
  A^2 \sqrt{N}.
\end{equation}
If
\begin{equation}
  \label{eq:7}
   \delta\geq \frac{c_0\log^q N}{w^{q}N^{q/2-1}}+2\exp\left(-c_1 a_0
 u^2\right)
\end{equation}
then $r_M(F,\tau, \delta)\leq A.$
\end{Corollary}
\proof
Clearly, it follows from Theorem~\ref{theo:multiplier_process_shahar} that if
\begin{equation*}
   \delta\geq \frac{c_0\log^q N}{w^{q}N^{q/2-1}}+2\exp\left(-c_1
 u^2\left(\frac{\ell^*(F\cap A D_{f^*})}{d_2(F\cap A D_{f^*})}\right)^2\right)
\end{equation*}
then $r_M(F,\tau, \delta)\leq A$. The claim follows because if $F \cap AD_{f^*}$ is nonempty,
\begin{equation*}
 \frac{\ell^*(F\cap A D_{f^*})}{d_2(F\cap A D_{f^*})} \geq a_0
\end{equation*}
for a suitable absolute constant $a_0$.
\endproof

\begin{Remark}
The estimate in Corollary \ref{coro:multiplier_process} can be rather loose. The reason for the suboptimal estimate is that usually, the Gaussian mean-width $\ell^*(F \cap A D_{f^*})$ is much larger than $d_2(F \cap A D_{f^*})$. For example, let $F=\{ \inr{t,\cdot} : t \in S^{d-1}\}$ be the class of linear functionals on $\R^d$ indexed by the Euclidean unit ball. Assume that $X$ is an isotropic vector -- that is, its covariance structure coincides with the standard Euclidean structure on $\R^d$; that $f^*=0$; and that $A \leq 1$. Then $F \cap AD_{f^*} = \{\inr{t,\cdot} : \|t\|_2 \leq A\}$, $d_2(F \cap AD_{f^*})=A$ and $\ell^*(F \cap AD_{f^*}) = A \sqrt{d}$, implying that
\begin{equation} \label{eq:Dvo-dim}
 \frac{\ell^*(F\cap A D_{f^*})}{d_2(F\cap A D_{f^*})} \geq \sqrt{d}
\end{equation}
which is significantly larger than an absolute constant.

Having said that, the question of an accurate probability estimate is not the main issue of this article and we will not explore that point further.
\end{Remark}

Next, we provide an estimate on $\gamma_\cO(\rho, \tau, \delta)$ that follows from Theorem~\ref{theo:multiplier_process_shahar} when $F$ is $L$-subgaussian and $\xi \in L_q$ for some $q>2$. The proof is identical to the one of Corollary~\ref{coro:multiplier_process}  and is omitted.

\begin{Corollary}
  \label{coro:lambda0}
Let $F$ be a closed, convex $L$-subgaussian class of functions and let
$\xi = Y-f^*(X)\in L_q$ for some $q>2$. Set $w, u>c$, $\tau>0$, $0<\delta<1$ and $\rho>0$. If $A>0$ satisfies
\begin{equation*}
  c_2 L w u \norm{\xi}_{L_q} \ell^*\left(F\cap K_\rho(f^*)\cap r(\rho)D_{f^*}\right)\leq \tau A \sqrt{N}.
\end{equation*}
and
\begin{equation*}
   \delta\geq \frac{c_0\log^q N}{w^{q}N^{q/2-1}}+2\exp\left(-c_1a_0 u^2\right)
\end{equation*}
then $\gamma_\cO(\rho, \tau, \delta)\leq A. $
\end{Corollary}

Finally, when $F$ is a $L$-subgaussian class it follows from a standard chaining argument (cf. \cite{Talagrand:05} or \cite{shahar_multi_pro}) that
\begin{equation}\label{eq:chaining_base}
\E \sup_{f\in F}\Big|\frac{1}{N}\sum_{i=1}^N \eps_i f(X_i)\Big|\leq \frac{c_0 L \ell^*(F)}{\sqrt{N}}.
\end{equation}
This observation will be used in what follows to upper bound $r_Q$.

\subsection{Estimates under a limited moment condition}
\label{sec:monotone-rearrangement}
In this section we shall consider the case of a class that need not be subgaussian, but rather, the growth of moments of class members is well-behaved up to some point. More accurately, we will assume that there is some $p_0$ for which, for every $f,h \in F \cup \{0\}$ and $2 \leq p \leq p_0$,
\begin{equation} \label{eq:local-limited-moments}
\|f-h\|_{L_p} \leq L \sqrt{p} \|f-h\|_{L_2}.
\end{equation}
In contrast, a subgaussian condition is equivalent to having  $\|f-h\|_{L_p} \leq L\sqrt{p}\|f-h\|_{L_2}$ for every $p \geq 2$.

The motivation for this type of limited moment assumption is the LASSO estimator. Recent results on properties of the {\it basis pursuit algorithm} in $\R^d$ \cite{LM_compressed,DLR} indicate that \eqref{eq:local-limited-moments} for $p_0 \sim \log d$ should suffice for an optimal estimate on the performance of the LASSO -- as if the class were subgaussian.

When analyzing the LASSO via the computation of the fixed points $r_M$ and $r_Q$, one encounters the following scenario. Let $X=(x_j)_{j=1}^d$ be a random vector in $\R^d$ and set $X_1,\ldots,X_N$ to be independent copies of $X$. Let $X_i(j)$ be the $j$-th coordinate of $X_i$ and thus $(X_{i}(j))_{i=1}^N$ is a random vector with independent coordinates, distributed as $x_j$.

Consider the random variables appearing in the definition of $r_M$ and $r_Q$ in the LASSO case:
\begin{equation} \label{eq:finite-max-1}
\max_{1 \leq i \leq d} \left|\sum_{i=1}^N \eps_i X_i(j)\right|,
\end{equation}
and
\begin{equation} \label{eq:finite-max-2}
\max_{1 \leq i \leq d} \left|\sum_{i=1}^N \eps_i \xi_i X_i(j)\right|
\end{equation}

The aim of this section is to derive upper bounds on \eqref{eq:finite-max-1} in expectation and \eqref{eq:finite-max-2} in deviation when each $x_j$ satisfies that
$$
\|x_j\|_{L_p} \leq L\sqrt{p}\|x_j\|_{L_2}
$$
for $p \lesssim \log d$. Note that an upper bound on the centered empirical process involved in the definition of $\gamma_\cO(\rho)$ will follow from a symmetrization argument and a bound on \eqref{eq:finite-max-2}.

The obvious difference between \eqref{eq:finite-max-1} and \eqref{eq:finite-max-2} are the multipliers
$(\xi_i)_{i=1}^N$: although the $x_j$'s have $\sim \log d$ moments, $\xi$ may be heavy-tailed, in the sense that it only belongs to $L_q$ for some fixed $q>2$; this difference makes the analysis of \eqref{eq:finite-max-2} more difficult.

Upper bounds on \eqref{eq:finite-max-1} and \eqref{eq:finite-max-2} are obtained under the following assumption.

\begin{Assumption} \label{ass:max-finite}
Let $N \leq d$, $t \geq 4$ and set $p_0=t \log d$. Assume that $p_0\lesssim N$ (and note that $p_0\geq \log N$) and that for every $1 \leq j \leq d$ and $p \leq p_0$, $\|x_j\|_{L_p} \leq L\sqrt{p}\|x_j\|_{L_2}$. Consider $\xi \in L_q$ for some $q>2$; let $r=\min\{1/2+q/4,2\}$; set $r^\prime$ to be the conjugate index of $r$; and assume that $4r^\prime\max\{2,1+a_0/a_1\} \leq t\log N$ (where $a_0$ and $a_1$ are two absolute constants to be specified later -- in Lemma~\ref{lemma:single-small-coordinates} and Lemma~\ref{lemma:large-coordinates}).
\end{Assumption}

Under this assumption we will prove the following:

\begin{Theorem}\label{thm:limited-moments}
Let the random vector $X$ and $\xi=Y-f^*(X)$ satisfy Assumption~\ref{ass:max-finite}.  Then,
\begin{equation}\label{eq:expected_quad}
\E \max_{1 \leq i \leq d} \left|\frac{1}{\sqrt{N}} \sum_{i=1}^N \eps_i X_i(j)\right| \leq c_0 \sqrt{\log d} \cdot L \max_{1 \leq j \leq d} \|x_j\|_{L_2}.
\end{equation}
Also, for every $u>2, v>0, w\geq2$ and for $p=p_0/2$ and $m =  p/\log(eN/p)$, one has that with probability at least
\begin{equation}\label{eq:final_result_multiproc_weak_moment}
1-\frac{\exp(-p/2)}{u^{2p}}-\frac{4 \exp(-p/2)}{u^{c_1 m}}-\frac{c_2 \log^q N }{w^q N^{q/2-1}}-2 \exp(-v^2t\log d),
\end{equation}
\begin{equation}\label{eq:upper_bound_deviation_lasso2}
\max_{1\leq j\leq d}\Big|\sum_{i=1}^N \eps_i \xi_i X_i(j)\Big| \leq c_3(q) (uw + u^2 v) L\|\xi\|_{L_q}  \sqrt{N}\sqrt{t\log d} \max_{1 \leq j \leq d} \|x_j\|_{L_2}.
\end{equation}
\end{Theorem}

The proofs of both estimates in Theorem~\ref{thm:limited-moments}  follow from a more general result, established in \cite{shahar_multi_pro}, on the supremum of a centered multiplier process under a limited moment assumption like \eqref{eq:local-limited-moments}. Although the estimate in \cite{shahar_multi_pro} is stated for the centered empirical process (cf. Section~4 there) its proof is actually based on an estimate on the symmetrized process. The proof of Theorem \ref{thm:limited-moments} will be presented in final section of this article.

\section{The LASSO under a limited moment assumption}
\label{sec:results-under-moment-lin}
In this section, we obtain complexity-dependent error rates for the LASSO. Our aim is to show that the LASSO (almost) achieves the minimax rates of convergence in the ``true model", and the meaning of the ``true model" is, in this case, the smallest $\ell_1^d$-ball centered in $0$ that contains  $t^*$. Thus, the price one has to pay for not knowing $\|t^*\|_1$ is rather minimal. 

The rate we shall be comparing the LASSO's performance to is the minimax rate of the following problem. Let $X\sim\cN(0,I_{d\times d})$ and set $\xi\sim\cN(0,\sigma^2)$ to be independent of $X$. Let $\rho>0$, consider an unknown $t_0 \in \rho B_1^d$ and  put $Y=\inr{X,t_0}+\xi$.

Let $c_0,c_1$ and $c_2$ be well-chosen absolute constants and  consider the cases $\log d\leq N\leq c_0 d$  or $c_1 d \leq N$. Following \cite{LM13}, if
\begin{equation*}
s_M^2(\rho) = c_2 \left\{
\begin{array}{cc}
\frac{\sigma^2 d}{N} & \mbox{ if } \rho^2 N \geq \sigma^2 d^2\\
\rho \sigma \sqrt{\frac{1}{N}\log \Big(\frac{e \sigma d}{\rho \sqrt{N}}\Big)} & \mbox{ if } \sigma^2 \log d \leq \rho^2 N \leq \sigma^2 d^2\\
\rho \sigma \sqrt{\frac{\log d}{N}} & \mbox{ if } \rho^2 N \leq \sigma^2 \log d
\end{array}
\right. \mbox{ and } s_Q^2(\rho)  \left\{
\begin{array}{cc}
 = 0 & \mbox{ if }  N \geq c_0 d\\
 \lesssim \rho^2/d & \mbox{ if } c_0 d\leq N \leq c_1 d\\
\sim \frac{\rho^2}{N}\log\Big(\frac{d}{N}\Big) & \mbox{ if } N\leq c_1 d,
\end{array}
\right.
\end{equation*}
then the minimax rate of convergence in the class $\rho B_1^d$ is 
\begin{equation}\label{eq:minimax_rate_ell_1}
\max\Big\{s_M^2(\rho), s_Q^2(\rho)\Big\}
\end{equation}
when $\rho \geq\sigma \sqrt{(\log d)/N}$ and $\rho^2$ when $\rho \leq\sigma \sqrt{(\log d)/N}$. Note that when $c_0d\leq N\leq c_1d$ (i.e. $N\sim d$), $s_Q^2(\rho)$ decays rapidly from $\frac{\rho^2}{N}\log({d}/{N})$ to $0$ and there are no precise estimates on the minimax rate in that range.

It turns out that for this problem -- the so-called Gaussian linear model -- the minimax rate in $\rho B_1^d$ is achieved by the Empirical Risk Minimization procedure (see, e.g., \cite{LM13}); however, an underlying assumption is that $\rho$ part of the information one is given. Thanks to regularization, and specifically, thanks to the LASSO, one does not need to know the value of $\norm{t_0}_1$ in advance to achieve the minimax rate, at least up to a logarithmic term. In fact, the optimal rate can be achieved using regularization in a much more general framework than just the Gaussian linear model -- as will be explained below. 

In what follows we will compare the rates obtained for the LASSO in Theorem~\ref{thm:main} in the high-dimensional case, that is, when $N\leq c_1 d$. One may do the same when $N \geq c_0d$ and we leave that to the reader.

Let $X$ be a random vector in $\R^d$ and consider the class of linear functionals
$F=\{f_t=\inr{\cdot,t}:t\in\R^d\}$. In particular, if $Y \in L_2$ is an arbitrary target random variable then
$f^*=f_{t^*}=\inr{\cdot,t^*}$ satisfies
\begin{equation}
  \label{eq:4}
  t^*\in\argmin_{t\in \R^d}\E\big(Y-\inr{X,t}\big)^2.
\end{equation}
As noted in the Introduction, the regularization function associated with the LASSO is the $\ell_1^d$-norm: for every $t=(t_j)_{j=1}^d \in\R^d$,
$$
\Psi(f_t)=\norm{t}_1=\sum_{j=1}^d |t_j|.
$$
Clearly, as a norm, the $\ell_1^d$-regularization function satisfies
Assumption~\ref{assum:reg-function} for $\eta=1$.

\vskip0.4cm

The LASSO with regularization parameter $\lambda$ produces
\begin{equation}
  \label{eq:lasso}
    \hat t \in \argmin_{t\in     \R^d}\Bigl(\frac{1}{N}\sum_{i=1}^N\big(Y_i-\inr{X_i,t}\big)^2+\lambda\norm{t}_1\Bigr),
\end{equation}
and one would like to control $\|f_{\hat{t}}-f_{t^*}\|_{L_2}^2=\E\inr{X,\hat t-t^*}^2$, where the expectation is taken with respect to $X$ conditionally to the data.

It should be noted that despite the LASSO's popularity, there are relatively few results in the random design scenario we are interested in (see, e.g., \cite{MR2981422}, \cite{MR2820635} and chapter~8.2 in \cite{MR2829871}). The overwhelming majority of existing results have been obtained for the linear model with subgaussian noise and a fixed design (i.e., each data point is of the form $Y_i=\inr{t^*,z_i}+\xi_i$) -- and the deterministic design matrix, whose rows are the vectors $z_i$, satisfies some form of the {\it Restricted Isometry Property} -- for example, the \textit{Restricted Eigenvalue Condition} (REC) from \cite{MR2533469} or the \textit{Compatibility Condition} (CC) from \cite{deter_lasso}).

To define the restricted eigenvalue condition, let us introduce the following notation:  for $x\in\R^d$ and a set $S_0\subset\{1,\ldots,d\}$ of cardinality  $|S_0|\leq s$, let $S_1$ be the set of indices of the $m$ largest coordinates of $(|x_i|)_{i=1}^d$ that are outside $S_0$. Let $x_{S_{01}}$ be the restriction of $x$ to the set $S_{01}=S_0\cup S_1$.

\begin{Definition}[\cite{MR2533469}]\label{def:rec}
  Let $\Gamma$ be an ${N\times d}$ matrix. For $c_0\geq1$ and an integer $1\leq s\leq m\leq d$ for which $m+s\leq d$, the \textbf{restricted eigenvalue constant} is
  \begin{equation*}
\kappa(s,m,c_0) =    \min\left\{\frac{\norm{\Gamma x}_2}{\norm{x_{S_{01}}}_2}:S_0\subset\{1,\ldots,d\}, |S_0|\leq s, \norm{x_{S_0^c}}_1\leq c_0 \norm{x_{S_0}}_1\right\}.
  \end{equation*}
The matrix $\Gamma$ satisfies the \textbf{Restricted Eigenvalue Condition (REC) of order $s$ with a constant $c$} if $\kappa(s,s,3)\geq c$.
\end{Definition}

One can show (see, \cite{MR2533469}, \cite{MR2807761}) that if $\Gamma$ satisfies REC and $\lambda\gtrsim \sigma \sqrt{(\log d)/N}$, then with high probability (with respect to the noise), simultaneously for every $1\leq p\leq 2$,
\begin{equation}
  \label{eq:lasso__rec_p}
  \norm{\hat t-t^*}_p^p\lesssim_p \norm{t^*}_0\left(\frac{\sigma}{\kappa(s,s,3)}\sqrt{\frac{\log d}{N}}\right)^p
\end{equation}where $\norm{t^*}_0$ is the cardinality of the support of $t^*$.

\vskip0.4cm

The main result in this section is an estimate on $\|\hat{f}-f^*\|_{L_2}^2$ that depends on $\|t^*\|_1$ rather than on the cardinality of the support of $t^*$ (we refer to \cite{LM_sparsity} for ``sparsity-dependent" rates of convergence for the LASSO in the same framework as we consider here). Such result follow from Theorem~\ref{thm:main}, and to that end, one has to construct a function $r(\cdot)$ as in \eqref{eq:r-rho} and to compute  $\lambda_0(\delta,\gamma)$ as in Definition~\ref{def:lambda0}. We will do so under the following situation:
Set $a_2 \geq 4$, $2\leq p_0=a_2 \log d\lesssim N$, $q>2$, $r=\min\{1/2+q/4,2\}$ and $r^\prime$ that is the conjugate index of $r$. Assume that $4r^\prime\max\{2,1+a_0/a_1\}\leq a_2 \log N$ (which is equivalent to assuming that $q>2+c_1/\log N$ for some constant $c_1=c_1(a_0,a_1, a_2)$).
Let $X=(x_j)_{j=1}^d$ be a random vector and note that the coordinates $x_1,...,x_d$ need not be independent.
\begin{Assumption} \label{ass:moments_design_noise}
Using the above notation, assume that there are constants $\kappa_0,\kappa$ and $\eps$ for which the following holds:
\begin{description}
\item{$\bullet$} For every $1\leq j\leq d$ and every $2\leq p\leq p_0$, $\norm{x_j}_{L_p}\leq \kappa_0 \sqrt{p}\norm{x_j}_{L_2}$.
\item{$\bullet$} $X$ satisfies a small-ball condition with constants $\kappa$ and $\eps$; that is, for every $t\in\R^d$,
 \begin{equation}
   \label{eq:6}
   Pr\left(|\inr{X,t}|\geq \kappa \norm{\inr{X,t}}_{L_2}\right)\geq \eps.
 \end{equation}
\item{$\bullet$} $\xi=Y-f^*(X) \in L_q$.
\end{description}
\end{Assumption}

To put this assumption in some perspective, note that an obvious underlying condition in any estimation problem with respect to the squared loss is that $\E(f(X)-Y)^2$ is defined for any $f\in F$, and in particular, that $\xi =Y-f^*(X) \in L_2$. Thus, assuming that $\xi \in L_q$ for some $q>2+c_1/\log N$ is not very restrictive. Also, as noted previously, the small-ball assumption is rather minimal.

The most restrictive component of Assumption~\ref{ass:moments_design_noise} is the moment assumption on the coordinates of $X$ -- that their moments exhibit a  subgaussian behavior, up to, roughly, $p \sim \log d$.

While this assumption can be weakened to other types of moment growth condition (e.g., $\norm{x_j}_{L_p}\leq \kappa_0 p^\alpha\norm{x_j}_{L_2}$ for some $\alpha\geq1/2$ and up to $p\sim\log d$), the resulting analysis is more involved (see \cite{LM_compressed}), and will not be explored here.

Finally, \cite{LM_compressed} shows that even if one assumes a subgaussian behavior of the coordinates $x_i$, but only up to $p \sim (\log d)/(\log \log d)$, Basis Pursuit may fail to recover even a $1$-sparse vector, implying that the choice of $p_0$ in Assumption \ref{ass:moments_design_noise} can not be relaxed significantly.

Given any $\rho\geq0$, set $M=\max_{1\leq j\leq d}\norm{x_j}_{L_2}$, let $\sigma_q = \norm{\xi}_q$ and put
\begin{equation*}
\Lambda(\rho) = \frac{\kappa_0\rho M}{\kappa^2 \eps} \sqrt{\frac{\log d}{N}}.
\end{equation*}
Moreover, for $R(t) = \E(Y-\inr{X,t})^2$, one has
\begin{equation*}
R(t) - R(t^*) = \E \inr{X,t-t^*}^2,
\end{equation*}
because $\inr{X,t^*}$ is the best approximation of $Y$ in a closed subspace of $L_2$. Thus, the estimation bounds also lead to excess risk bounds.

\begin{Theorem}  \label{theo:lasso-small-ball-weak-moments}
There are absolute constants $c_0,...,c_{6}$ for which the following holds.
  Assume that $X$ and $\xi=Y-f^*(X)$ satisfy
  Assumption~\ref{ass:moments_design_noise} and that $N \leq d$.  Let $u>2, v>0$ and $w\geq 2$, and set $p=(a_2/2) \log d$ and  $m = p/\log(eN/p)$. Put
  \begin{equation}
    \label{eq:3}
    \delta= \frac{\exp(-p/2)}{u^{2p}}-\frac{4 \exp(-p/2)}{u^{c_0 m}}-\frac{c_1 \log^q N }{w^q N^{q/2-1}}-2 \exp(-v^2t\log d)
  \end{equation}
and set
\begin{equation*}
r^2(\rho)=c_2\left\{
\begin{array}{cc}
(uw+u^2w) \sigma_q \Lambda(\rho) & \mbox{ if } N\geq (\kappa\eps/32)^2 d\\
\max\Big\{(uw+u^2v)\sigma_q \Lambda(\rho), \kappa^2\Lambda^2(\rho)\Big\}&\mbox{ otherwise.}
\end{array}
\right.
\end{equation*}
If $\hat t$ is produced by the LASSO for a regularization parameter
  \begin{equation*}
    \lambda > c_4(uw+u^2v)
   \kappa_0 \norm{\xi}_{L_q} \eta^3 M\sqrt{\frac{\log d}{N}},
  \end{equation*}
  then with probability at least  $1-5\delta-2\exp(-\eps^2N/2)$,
\begin{equation*}
  R(\hat t)-R(t^*) = \norm{\inr{X,\hat t-t^*}}_{L_2}^2\leq c_5
  \max\Big\{r^2(c_{6} \norm{t^*}_1), \frac{\lambda}{\kappa^2\eps} \norm{t^*}_1\Big\}.
\end{equation*}
\end{Theorem}

Observe that like known estimates on the LASSO, and despite imposing considerably weaker assumptions on $X$ and $Y$, the regularization parameter in Theorem~\ref{theo:lasso-small-ball-weak-moments} is of the order of $\|\xi\|_{L_q}\sqrt{(\log d)/N}$. And, when $\|\xi\|_{L_q}$ is equivalent to $\sigma$ -- the variance of $\xi$ -- then for $N\lesssim d$, the rate of convergence is
\begin{equation*}
 c(M) \max\left\{\sigma \norm{t^*}_1 \sqrt{\frac{\log d}{N}}, \norm{t^*}_1^2\frac{\log d}{N}\right\}
\end{equation*}
for a constant that depend only on $M$. 

Hence, up to a logarithmic factor, the LASSO attains the minimax rate in $\norm{t^*}_1B_1^d$ when $\log d \leq N\lesssim d$ and when $\norm{t^*}_1\geq \sigma \sqrt{\log d/N}$; moreover, it does so without knowing in advance the identity of the  ``true model" $\norm{t^*}_1 B_1^d$.

Note that one may want to combine the sparsity-dependent error rate from Theorem~1.3 in \cite{LM_sparsity} and the complexity-dependent error rate from Theorem~\ref{theo:lasso-small-ball-weak-moments}. To simplify the exposition, results from \cite{LM_sparsity}  have been stated under a subgaussian assumption on the design. Therefore, we will also make this assumption below. Note also that the probability estimate from Theorem~\ref{theo:lasso-small-ball-weak-moments} can be improved under the subgaussian assumption on the design and that the third condition from Assumption~\ref{ass:moments_design_noise} (i.e. $q>2 + c_1/\log N$) can be relaxed to only $q>2$ (see more details in the next section, and, in particular,  Theorem~\ref{theo:ell_p_reg}). Combining the two approaches, one has that when $X$ is isotropic and $L$-subgaussian, and when $\xi\in L_q$ for some $q>2$ then for any $u,w>c$ with probability larger than $1-\delta$ for
 \begin{equation*}
 \delta = 2\exp(- c_2 N/L^8)- \frac{c_0\log^q N}{w^q N^{q/2-1}}-c_0 \exp(-c_1u^2/L^2),
 \end{equation*}
 the LASSO estimator  $\hat t$ with the universal regularization parameter $\|\xi\|_{L_q} \sqrt{(\log d)/N}$ satisfies that
\begin{equation}\label{eq:combined_sparse_comp}
\norm{\hat t - t^*}_2^2\lesssim_{L, q} \min\left\{\frac{\norm{t^*}_0 \sigma^2 \log d}{N}, \max\left\{\sigma \norm{t^*}_1 \sqrt{\frac{\log d}{N}}, \norm{t^*}_1^2\frac{\log d}{N}\right\} \right\}
\end{equation}
when $N \gtrsim \norm{t^*}_0 \log(d/\norm{t^*}_0)$.

Note that seemingly, \eqref{eq:combined_sparse_comp} exhibits a different rate than  Corollary~9.1 in \cite{MR2829871} (see also \eqref{eq:vlad_lasso} above): the extra (and necessary) $\norm{t^*}_1^2\frac{\log d}{N}$ term in \eqref{eq:combined_sparse_comp}, which appears only in the random design scenario. As a result, the rates of convergence of the LASSO appears to deteriorate when
\begin{equation}\label{eq:regime_diff_rand_deter}
\sigma\sqrt{\frac{N}{\log d}}\leq \norm{t^*}_1\leq \sigma \sqrt{\norm{t^*}_0}.
\end{equation}
However, the sparsity-dependent error rate, and therefore Equation~\eqref{eq:combined_sparse_comp},  holds only when $N \gtrsim \norm{t^*}_0 \log(d/\norm{t^*}_0)$. And, when $N \gtrsim \norm{t^*}_0 \log d$ (which is only slightly larger than $\norm{t^*}_0 \log (d/\norm{t^*}_0)$), the error rates in the two scenarii (random and deterministic design) are the same and are given by
\begin{equation*}
\min\left\{\frac{\sigma^2 \norm{t^*}_0 \log d}{N}, \sigma \norm{t^*}_1 \sqrt{\frac{\log d}{N}} \right\}.
\end{equation*}

\vspace{0.8cm}

\noindent {\bf Proof of Theorem \ref{theo:lasso-small-ball-weak-moments}.}
 As noted previously, since $\norm{\cdot}_1$ is a norm, $\Psi(t)=\|t\|_1$ satisfies Assumption~\ref{assum:reg-function} for $\eta=1$; therefore, Theorem~\ref{thm:main} may be applied here, and one has to control $r(\rho) \equiv \max\{r_M(\rho), r_Q(\rho)\}$ and $\lambda_0(\delta,\tau)$. In what follows we will invoke the results of Section~\ref{sec:pre:estimates} and estimate these parameters.

Set $F(f^*,\rho)=  F\cap K_{\rho}(f^*)-f^*$ and recall that $r_Q(\rho)=r_Q(F\cap K_\rho(f^*), \kappa\eps/32)$ is determined by the behavior of
\begin{equation}
  \label{eq:complexity-term-moments}
  (\star)= \E\sup_{f\in F(f^*,\rho)\cap r D}\left|\frac{1}{\sqrt{N}}\sum_{i=1}^N \eps_i
  f(X_i)\right|;
\end{equation}
as a consequence, it suffices to upper bound $(\star)$. Let $\cE=\{t\in\R^d:\E\inr{X,t}^2\leq1\}$, put $\cE^\circ$ to be the polar of $\cE$ (that is, $\cE^\circ = \{ u: \sup_{t \in \cE} |\inr{u,t}| \leq 1\}$), and set $\norm{t}_\cE=\sup_{x\in\cE}\inr{x,t}$. Thus,
\begin{align*}
  &(\star)=\E\sup_{t\in \rho B_1^d\cap r
    \cE}\left|\frac{1}{\sqrt{N}}\sum_{i=1}^N \eps_i \inr{X_i,t}\right| \leq \min \left\{ \E \sup_{t\in \rho B_1^d}\left|\frac{1}{\sqrt{N}}\sum_{i=1}^N \eps_i \inr{X_i,t}\right|, \E \sup_{t\in r \cE}\left|\frac{1}{\sqrt{N}}\sum_{i=1}^N \eps_i \inr{X_i,t}\right| \right\}
\\
&= \min\left\{\rho \E\norm{\frac{1}{\sqrt{N}}\sum_{i=1}^N \eps_i X_i}_{\ell_\infty^d},r\E\norm{\frac{1}{\sqrt{N}}\sum_{i=1}^N \eps_i X_i}_{\cE^\circ}\right\}.
\end{align*}

It is standard to verify (see, for instance, the proof of Lemma~2.2 in \cite{LM_lin_agg}) that
\begin{align*}
  \E\norm{\frac{1}{\sqrt{N}}\sum_{i=1}^N \eps_i X_i}_{\cE^\circ}\leq \sqrt{d}.
\end{align*}
Moreover, by (\ref{eq:expected_quad}),
\begin{align*}
  \E\norm{\frac{1}{\sqrt{N}}\sum_{i=1}^N \eps_i
    X_i}_{\ell_\infty^d} = & \E \max_{1 \leq j \leq d} \left|\frac{1}{\sqrt{N}} \sum_{i=1}^N \eps_i X_i(j)\right|
    \leq  c_0 \kappa_0 \sqrt{\log d}\max_{1\leq
    j\leq d}\norm{x_j}_{L_2}.
\end{align*}
Therefore,
\begin{equation}
  \label{eq:4'}
  (\star)\leq \min\left\{c_0\rho \kappa_0 \sqrt{\log d}\max_{1\leq
    j\leq d}\norm{x_j}_{L_2}, r\sqrt{d} \right\},
\end{equation}
and setting $\gamma=\kappa\eps/32$, one has
\begin{equation*}
  r_Q(\rho) \leq \left\{
    \begin{array}{cc}
      0 & \mbox{ if } N\geq \gamma^2 d\\
       \frac{c_0 \rho\kappa_0}{\gamma} \sqrt{\frac{\log d}{N}}M & \mbox{ otherwise.}
    \end{array}
\right.
\end{equation*}

Next, let us establish a high probability upper bound on $r_M(\rho) = r_M(F\cap K_\rho(f^*),
\kappa^2\eps/80, \delta/4)$. Note that
\begin{align*}
  \phi_N(F\cap K_\rho(f^*), f^*,s) = & \sup_{t\in \rho B_1^d\cap s
    \cE}\left|\frac{1}{\sqrt{N}}\sum_{i=1}^N \eps_i \xi_i
  \inr{X_i,t}\right|
  \leq  \rho \max_{1\leq j\leq
    d}\left|\frac{1}{\sqrt{N}}\sum_{i=1}^N\eps_i \xi_i X_i(j)\right|.
\end{align*}
Applying the second result of Theorem~\ref{thm:limited-moments} for $u>2, v>0$, $w\geq2$, $p=(a_2/2) \log d, m = p/\log(eN/p)$ and
\begin{equation}
  \label{eq:13}
  \delta=\frac{\exp(-p/2)}{u^{2p}}-\frac{4 \exp(-p/2)}{u^{c_0 m}}-\frac{c_1 \log^q N }{w^q N^{q/2-1}}-2 \exp(-v^2t\log d),
\end{equation}
it follows that with probability at least $1-\delta$,
\begin{equation*}
   \phi_N(F\cap K_\rho(f^*), f^*,s)\leq c_2 \kappa_0 (uw+u^2 v) \norm{\xi}_{L_q}\rho  M
  \sqrt{\log d};
\end{equation*}
thus,
\begin{equation*}
  r_M^2(\rho)\leq \frac{c_3\kappa_0 (uw+u^2 v)}{\kappa^2 \eps} \norm{\xi}_{L_q}\rho M
  \sqrt{\frac{\log d}{N}}.
\end{equation*}
\vskip0.5cm

Finally, let us identify an upper bound on $\lambda_0(\delta,\tau)$ for $\tau = 3/(80 \eta^3)$.  Let $\{e_1,\ldots,e_d\}$ be the canonical basis of $\R^d$. Since $K_\rho(f^*)=\{t : \|t-t^*\|_1 \leq \rho\}$, we have
\begin{align*}
 &(\star_1) = \sup_{f\in F\cap
K_{\rho}(f^*)\cap r(\rho)D_{f^*}} \left(\frac{1}{N}\sum_{i=1}^N
\xi_i(f-f^*)(X_i) - \E \xi (f-f^*)(X)\right)\\
&\leq \rho  \max_{t-t^*\in\{\pm e_1,\ldots,\pm
   e_d\}}\left(\frac{1}{N}\sum_{i=1}^N \xi_i
 \inr{X_i,t-t^*}-\E\xi\inr{X,t-t^*}\right)\\ & = \rho  \max_{t\in\{\pm e_1,\ldots,\pm
   e_d\}}\left(\frac{1}{N}\sum_{i=1}^N \xi_i
 \inr{X_i,t}-\E\xi\inr{X,t}\right).
\end{align*}
Recall that $X=(x_j)_{j=1}^d$. By a standard symmetrization argument
(see, for example, Lemma~2.3.7 in \cite{vanderVaartWellner}), if $z\geq 4
\max_{1\leq j\leq d}\sqrt{{\rm Var}(\xi x_j)/N}$ then
\begin{align*}
  Pr \Big(\max_{t\in\{\pm e_1,\ldots,\pm e_d\}}&\frac{1}{N}\sum_{i=1}^N
  \xi_i \inr{X_i,t}-\E\xi\inr{X,t}\geq z\Big)
  \leq 4 Pr\Big(\max_{t\in\{\pm
    e_1,\ldots,\pm e_d\}}\frac{1}{N}\sum_{i=1}^N \eps_i \xi_i
  \inr{X_i,t}\geq \frac{z}{4}\Big).
\end{align*}
Note that $\sqrt{{\rm Var}(\xi x_j)}\leq \sqrt{\E \xi^2 x_j^2}\leq \norm{\xi}_{L_q}\norm{x_j}_{L_{2q^\prime}}$ where $q^\prime$ is the conjugate index of $q/2$. Therefore, $\sqrt{{\rm Var}(\xi x_j)}\leq \kappa_0 \sqrt{2q^\prime} \norm{\xi}_{L_q}M$ as long as $2q^\prime\leq a_2 \log d$, i.e., when $q\geq 2+2/(a_2\log d-1)$ -- which is the case under Assumption~\ref{ass:moments_design_noise}. Therefore, applying the second result of Theorem~\ref{thm:limited-moments} for $\delta$ as in (\ref{eq:13}), it follows that with probability at least $1-\delta$,
\begin{equation*}
  (\star_1) \leq \rho c_0 \kappa_0 (uw+u^2v) \norm{\xi}_{L_q} M
  \sqrt{\frac{\log d}{N}},
\end{equation*}
and, for $\tau = 3/(80 \eta^3)$ one may select
\begin{equation*}
  \lambda_0(\delta,\tau)=c_4\kappa_0 (uw+u^2v) \norm{\xi}_{L_q} \eta^3 w M
  \sqrt{\frac{\log d}{N}}.
\end{equation*}
\endproof

\section{Regularization methods for subgaussian classes}
\label{sec:examples}
In this section we assume that $X$ is a random vector that takes its values in a Hilbert space $\cH$. The main examples we shall consider are when $\cH$ is the
$d$-dimensional Euclidean space and when it is the space of $m\times T$ matrices
endowed with the Frobenius norm.

The inner product in $\cH$ is denoted by $\inr{\cdot,\cdot}$, and the norm and unit ball endowed by the inner product are denoted by  $\norm{\cdot}_{\cH}$ and $B_{\cH}=\{t\in \cH:
\norm{t}_{\cH}\leq 1\}$ respectively.

There is another natural Hilbertian structure on $\cH$, endowed by $\Sigma=\E X X^\top$, the covariance operator associated with the random vector $X$. The corresponding unit ball is $\cE=\{t\in\cH:\E\inr{X,t}^2\leq1\}$, is an ellipsoid in $\cH$.

\vskip0.4cm

Let $T\subset\cH$ be a closed and convex set and put
\begin{equation*}
  t^*\in\argmin_{t\in T}\E(Y-\inr{X,t})^2;
\end{equation*}
thus, $\inr{X,t^*}$ is the best $L_2(\mu)$-approximation of $Y$ by a
linear functional $\inr{t,\cdot}$ for $t \in T$.

Let $\Psi(\cdot)$ be a regularization function on $\cH$ that satisfies
Assumption~\ref{assum:reg-function}. The goal is to estimate $t^*$  in $L_2(\mu)$ with a rate depending on $\Psi(f^*)$. To that end, set
\begin{equation}
  \label{eq:2}
  \hat t\in\argmin_{t\in T}\Big(\frac{1}{N}\sum_{i=1}^N\big(Y_i-\inr{X_i,t}\big)^2+\lambda\Psi(t)\Big)
\end{equation}
for a well-chosen regularization parameter $\lambda$.

Unlike the results of the previous section, in what follows we will assume that $F= \{\inr{t,\cdot} : t \in T\}$ is an $L$-subgaussian class (see Definition~\ref{def:subgaussian-class}). Moreover, $F$ satisfies a small-ball property with constants that depend only on $L$. Indeed, observe that for every $t \in T$
$$
\norm{\inr{X,t}}_{L_4}\lesssim L \norm{\inr{X,t}}_{L_2},
$$
and applying the Paley-Zygmund inequality (see, e.g., Corollary~3.3.2 in \cite{MR1666908}),
\begin{equation}
  \label{eq:small-ball-lin}
  Pr\left(|\inr{X,t}|\geq \kappa \norm{\inr{X,t}}_{L_2}\right)\geq \eps \mbox{ for } \kappa = 1/2 \mbox{ and } \eps = c/L^4.
\end{equation}

From here on we will say that the random vector $X$ taking its values in $\cH$ is $L$-subgaussian if the class consisting of all the linear functionals on $\cH$, i.e., $\{\inr{t,\cdot}: t \in \cH\}$, is $L$-subgaussian. Also, throughout this section, we will assume that $\xi=Y-f^*(X) \in L_q$ for some $q>2$, $\sigma_q = \|\xi\|_{L_q}$, and
$$
T\cap K_\rho(t^*)=\{ t \in T : \Psi(t-t^*) \leq \rho\}.
$$

\subsection{`Heavy tailed' noise}
\label{sec:results-sub-gaussian-lin}
Thanks to the subgaussian assumption, both $r(\rho)$ and $\lambda_0=\lambda_0(\delta, 3/(80\eta^3))$ may be determined using the Gaussian mean-widths of the sets $T \cap K_\rho(t^*)$ for all $\rho>0$. Recall that for $T_0\subset \cH$ the Gaussian mean-width of $T_0$ is $\ell^*(T_0) = \E \sup_{t\in T_0} G_t$, where $(G_t)_{t\in T_0}$ is the centered canonical Gaussian process indexed by $T_0$ with covariance structure given by $\E G_{t_1} G_{t_2} = \E \inr{X, t_1}\inr{X, t_2}$ for every $t_1,t_2\in T$.

\begin{Definition} \label{def:reg_params}
Let $r\cE_{t^*}=\left\{t\in\cH: \|\inr{t-t^*,\cdot}\|_{L_2(\mu)} \leq
  r\right\}=t^*+r\cE$, and for $\alpha,\beta>0$ set
  \begin{equation*}
    \tilde{r}_Q(\rho, \alpha) = \inf\left\{r>0:\ell^*\left(T\cap K_\rho(t^*)\cap r \cE_{t^*}\right)\leq
    \alpha r \sqrt{N}\right\}
  \end{equation*}
  and
  \begin{equation*}
    \tilde{r}_M(\rho, \beta) = \inf\left\{r>0:\ell^*\left(T\cap K_\rho(t^*)\cap r \cE_{t^*}\right)\leq
    \beta r^2 \sqrt{N}\right\}.
  \end{equation*}
\end{Definition}

Let $c_0$ be an absolute constant to be specified later. Fix $u,w>c$, and $\eps$ and $\kappa$ as in \eqref{eq:small-ball-lin}. Consider
\begin{equation}
  \label{eq:8}
  \alpha=\frac{\kappa\eps}{c_0L}, \ \  \beta = \frac{\kappa^2
    c_1\eps}{ L wu
    \norm{\xi}_{L_q}}, \ \ {\rm and } \ \ \gamma=c_0 \eta^3 L wu\norm{\xi}_{L_q},
\end{equation}
put
\begin{equation}\label{eq:r_sub_gaussian}
r(\rho)\geq \max\left\{\tilde{r}_Q(\rho, \alpha), \tilde{r}_M(\rho, \beta)\right\}
\end{equation}
and set
\begin{equation}
  \label{eq:12}
  \lambda_0(\gamma)=\gamma \sup_{\rho>0, t^*\in T}\frac{\ell^*(T\cap
    K_\rho(t^*)\cap r(\rho)\cE_{t^*})}{\rho \sqrt{N}}.
\end{equation}
The first result we present is rather general and holds for any closed and convex subset $T\subset\cH$ and any regularization function satisfying Assumption~\ref{assum:reg-function}. It allows to take into account an additional constraint on the ``signal'' $t^*\in T$.

\begin{Theorem}\label{theo:linear-regression}
There are absolute constants $c,c_1$ and $c_2$ for which the following holds. Let $\Psi$ be a regularization function satisfying Assumption~\ref{assum:reg-function}. Assume that $X$ is $L$-subgaussian for some $L>0$ and that $\xi=Y-\inr{X,t^*}$ is in $L_q$ for some $q>2$.

If $\hat t$  is given by \eqref{eq:2} for a regularization parameter $\lambda> \lambda_0(\gamma)$ as in \eqref{eq:12}, then with probability larger than
 \begin{equation}\label{eq:proba1}
 1-2\exp(- N \eps^2/8)- \frac{c_0\log^q N}{w^q N^{q/2-1}}-c_0 \exp(-c_1u^2/L^2),
 \end{equation}
\begin{equation*}
\norm{\inr{X,\hat t-t^*}}_{L_2}^2\leq \max\big\{
r(10 \eta\Psi(t^*))^2,(32/\kappa^2\eps)\lambda \Psi(t^*)\big\}
\end{equation*}
for $r(\cdot)$ given by \eqref{eq:r_sub_gaussian}.
\end{Theorem}

\proof The proof follows from Theorem~\ref{thm:main} by estimating $r(\rho)$ and $\lambda_0$ using the `local' Gaussian mean-widths of the sets $T\cap K_\rho(t^*)$.

Since $X$ is $L$-subgaussian, the process $\left\{\inr{X,t} : t \in \cH\right\}$ is $L$-subgaussian. Setting $F=\{\inr{t,\cdot} : t \in \cH\}$ and $f^*=\inr{t^*,\cdot}$, a standard chaining argument shows that
\begin{equation*}
  \E\sup_{f\in F\cap K_\rho(f^*)\cap r D_{f^*}}\Big|\frac{1}{N}\sum_{i=1}^N \eps_i
  (f-f^*)(X_i)\Big|\leq c_0L \frac{\ell^*(T\cap K_\rho(t^*)\cap r \cE_{t^*})}{\sqrt{N}}.
\end{equation*}
Thus,
\begin{equation}
  \label{eq:9}
  r_Q\Big(F\cap K_\rho(f^*),\frac{\kappa\eps}{32}\Big)\leq \tilde{r}_Q(\rho, \alpha).
\end{equation}

As for the fixed point associated with the multiplier process, it
follows from Corollary~\ref{coro:multiplier_process} that
\begin{equation}
  \label{eq:10}
  r_M\Big(F\cap K_\rho(f^*),\frac{\kappa^2\eps}{160},
  \frac{\delta}{4}\Big)\leq \tilde{r}_M(\rho,\beta)
\end{equation}for $\beta$ as defined in \eqref{eq:8}, and as long as
\begin{equation*}
  \frac{\delta}{4}\geq \frac{c_0\log^q N}{w^q N^{q/2-1}}+2 \exp(-c_1u^2/L^2).
\end{equation*}

Finally by Corollary~\ref{coro:lambda0}, $\lambda_0(\delta,\gamma)\leq \lambda_0(\gamma)$. The claim now follows from Theorem~\ref{thm:main}.
\endproof

If one is to apply Theorem \ref{theo:linear-regression}, an essential component is an upper bound on $\ell^*(T\cap K_\rho(t^*)\cap r \cE_{t^*})$ -- which in turn determines $r$ and $\lambda$. To simplify the analysis we shall use an additional assumption on $\Psi$:

\begin{Assumption} \label{ass:additional-on-psi}
Assume that for every $x,y\in\cH$ and $\lambda\geq0$,
\begin{equation}
  \label{eq:condition-norm}
 \Psi(x)=\Psi(-x), \hspace{0.5cm} \Psi(x+y)\leq \eta\big(\Psi(x)+\Psi(y)\big) \mbox{ and }
  \Psi(\lambda x)\leq \lambda \Psi(x).
\end{equation}
\end{Assumption}
Also, recall that $\cE=\{t\in\cH:\E\inr{X,t}^2\leq1\}$, $\sigma_q= \norm{\xi}_q$ and set $K=\{t\in\cH:\Psi(t)\leq1\}$.

\begin{Theorem} \label{theo:learning-linear-norm}
Assume that $\Psi$ satisfies Assumption \ref{ass:additional-on-psi} and that the assumptions of  Theorem~\ref{theo:linear-regression} hold. Let $\Lambda(\rho)\geq \rho \ell^*(K)/\sqrt{N}$ for every $\rho>0$, $w,u>c$ and consider the RERM
\begin{equation*}
  \hat t \in\argmin_{t\in \cH}\Big(\frac{1}{N}\sum_{i=1}^N(Y_i-\inr{X_i,t})^2+c_0\eta^3Lwu\sigma_q\Lambda(\Psi(t))\Big).
\end{equation*}
Then, with probability larger than the one in \eqref{eq:proba1},
\begin{equation*}
  R(\hat t)-R(t^*)=\norm{\inr{X,\hat t-t^*}}_{L_2}^2\leq c_0 r^2(10\eta\Psi(t^*))
\end{equation*}where, for $\alpha$ and $\beta$ defined in \eqref{eq:8} and, for any $\rho\geq0$,
\begin{equation}
  \label{eq:choice-of-s} r^2(\rho)=\left\{
    \begin{array}{cc}
\frac{ \Lambda(\rho)}{\beta} & \mbox{ if } N\geq \big(\ell^*(\cE)/\alpha\big)^2\\
\max\Big\{\frac{\Lambda(\rho)}{\beta} ,\frac{\Lambda^2(\rho)}{\alpha^2}\Big\}& \mbox{ otherwise.}
    \end{array}
\right.
\end{equation}
\end{Theorem}
\proof The result follows immediately from Theorem~\ref{theo:linear-regression}. Indeed, for every $\rho>0$, $r>0$ and $t^*\in T=\cH$,
\begin{equation*}
  \ell^*(T\cap K_\rho(t^*)\cap r \cE_{t^*})=\ell^*(K_\rho(0)\cap r
  \cE)\leq \ell^*(\rho K\cap r \cE)\leq \min\left\{\rho \ell^*(K), r \ell^*(\cE)\right\},
\end{equation*}
because $K_\rho(0)\subset \rho K=\{\rho t:t\in K\}$.

\endproof

Note that in a $d$-dimensional space, the trivial bound $\ell^*(\cE)\leq \sqrt{d}$ holds (see, e.g., Lemma~2.2 in \cite{LM_lin_agg}). Therefore,  one only needs to control $\ell^*(K)$. In the next section, we provide several examples of applications of Theorem~\ref{theo:learning-linear-norm} that follow from estimates on $\ell^*(K)$. We will simplify the analysis by assuming that there is some compatibility between the norm  $\norm{\cdot}_\cH$ and the one endowed by the covariance structure of $X$:
\begin{Assumption}\label{eq:isotropicity}
Assume that $X$ is isotropic; that is, for every $t\in\cH$, $\big(\E \inr{X,t}^2\big)^{1/2}=\norm{t}_{\cH}$.
\end{Assumption}
Observe that under Assumption~\ref{eq:isotropicity}, $\ell^*(K) = \E \sup_{t\in K} G_t$, where $(G_t)_{t\in K}$ is the canonical Gaussian process indexed by $K$ with the covariance $\E G_{t_1} G_{t_2} = \inr{t_1,t_2}$ for every $t_1,t_2\in K$, because the inner-product in $\cH$ coincides with the one endowed by $L_2(\mu)$.

\subsection{Regularization methods in  $\R^d$}
Consider a regularization function $\Psi(\cdot)$ satisfying Assumption~\ref{ass:additional-on-psi}. Assume that $X$ is $L$-subgaussian and isotropic in $\R^d$ with respect to the standard Euclidean inner-product, and that $\xi\in L_q$ for some $q>2$. Let $u,w>c$. For any $\rho\geq0$ set $\Lambda(\rho)\geq \rho \ell^*(K)/\sqrt{N}$ and put
\begin{equation*}
r^2(\rho)\sim_{L, q} \left\{
\begin{array}{cc}
 wu \sigma_q \Lambda(\rho) & \mbox{ when } N\gtrsim_L d
 \\
 \\
  \max\left\{wu \sigma_q \Lambda(\rho),\Lambda^2(\rho)\right\} & \mbox{ otherwise.}
\end{array}\right.
\end{equation*}
It follows from Theorem~\ref{theo:learning-linear-norm} that if
\begin{equation}\label{eq:reg-meth}
  \hat t \in\argmin_{t\in
    \R^d}\Big(\frac{1}{N}\sum_{i=1}^N(Y_i-\inr{X_i,t})^2+c_0\eta^3 L wu \sigma_q \Lambda(\Psi(t))\Big)
\end{equation}
then with probability larger than the one in \eqref{eq:proba1}
\begin{equation*}
  \norm{\inr{\hat t-t^*,\cdot}}_{L_2(\mu)}^2\lesssim r^2(10\eta \Psi(t^*)).
\end{equation*}
As a consequence, one can derive an estimation result for \eqref{eq:reg-meth} whenever $\ell^*(K)$ may be controlled from above. In the following section, we shall apply this observation to some classical problems and compare the error rates obtained by the RERM \eqref{eq:reg-meth} to the minimax rate in the ``true model" $\{t\in T: \Psi(t)\leq \Psi(t^*)\}$. 

\vskip0.4cm
\noindent{\bf Example: $\ell_p$-regularization for  $1\leq p\leq \infty$.}
In this section, we consider a regularization function $\Psi(t)=\norm{t}_p$ for some $p\geq1$. Assumption~\ref{ass:additional-on-psi} holds with $\eta=1$ because $\norm{\cdot}_p$ is a norm. In order to apply the general result for the RERM in \eqref{eq:reg-meth}, one has to compute the Gaussian mean-width of the unit ball associated with the regularization function $\Psi(\cdot)=\norm{\cdot}_p$.

In the range $1\leq p\leq 1+(\log d)^{-1}$, we recover the same result as for the LASSO, because $B_1^d\subset B_p^d\subset c B_1^d$ for a suitable absolute constant $c$; hence, $\ell^*(B_p^d)\sim \ell^*(B_1^d )\sim \sqrt{\log(ed)}$.

When $1+(\log(ed))^{-1}\leq p$, set $r$ to be the conjugate index for $p$ and one may easily verify that $\ell^*( B_p^d)\sim \sqrt{r} d^{1/r}$.

Applying Theorem~\ref{theo:learning-linear-norm}, one has the following:

\begin{Theorem}\label{theo:ell_p_reg}
Under the assumptions of Theorem \ref{theo:learning-linear-norm} and using its notation,
  \begin{itemize}
  \item
If $1\leq p\leq 1+1/(\log d)$ and
  \begin{equation*}
      \hat t \in\argmin_{t\in
    \R^d}\Big(\frac{1}{N}\sum_{i=1}^N(Y_i-\inr{X_i,t})^2+c_2 \eta_p^3 L wu\sigma_q\norm{t}_p\sqrt{\frac{\log d}{N}}\Big)
  \end{equation*}
then with probability larger than the one in \eqref{eq:proba1},
  \begin{equation*}
    \norm{\hat t-t^*}_2^2\lesssim_{p,L,q}\left\{
      \begin{array}{cc}
        wu\sigma_q\norm{t^*}_p\sqrt{\frac{\log d}{N}} &    \mbox{ if } N
        \gtrsim_L d,
        \\
        \\
        \max\left\{wu\sigma_q\norm{t^*}_p\sqrt{\frac{\log d}{N}},
        \norm{t^*}_p^2\frac{\log d}{N}\right\} & \mbox{ otherwise}.
      \end{array}
\right.
  \end{equation*}
\item If $p\geq 1+1/(\log d)$ and
  \begin{equation*}
      \hat t \in\argmin_{t\in
    \R^d}\Big(\frac{1}{N}\sum_{i=1}^N(Y_i-\inr{X_i,t})^2+c_2 \sigma_q L wu\norm{t}_p\frac{\sqrt{p/(p-1)}d^{(p-1)/p}}{\sqrt{N}}\Big),
  \end{equation*}
then with probability larger than the one in \eqref{eq:proba1}
  \begin{equation*}
    \norm{\hat t-t^*}_2^2\lesssim_{L,q} \left\{
      \begin{array}{cc}
        wu\sigma_q\norm{t^*}_p\frac{d^{(p-1)/p}}{p\sqrt{N}} &    \mbox{ if } N
        \gtrsim_L d,
        \\
        \\
        \max\left\{wu\sigma_q\norm{t^*}_p\frac{d^{(p-1)/p}}{\sqrt{N}},
        \norm{t^*}_p^2\frac{d^{2(p-1)/p}}{N}\right\} & \mbox{ otherwise}.
      \end{array}
\right.
  \end{equation*}
  \end{itemize}
\end{Theorem}

\begin{remark}[\textbf{the case $0<p<1$}]
Despite being a non-convex function, $\ell_p$-regularization for $0<p < 1$ has attracted much attention in the context of Signal Processing and High-Dimensional Statistics. Among the problems studied using $\ell_p$ regularization were the linear regression model with a deterministic design (cf. \cite{MR2882274,MR2816337,MR2879672}); the sequence space model \cite{MR1278886, MR2281879}; and the random design linear regression model \cite{WPGaoYang12}.

From our point of view, there is no particular restriction on $p$ as long as the regularization function satisfies Assumption~\ref{ass:additional-on-psi}. We can therefore consider a regularization function $\Psi(t)=\norm{t}_p$ for any $0<p<1$. In that range of $p$, Assumption~\ref{ass:additional-on-psi} holds for $\eta=\eta_p=2^{1/p}$ (see, for example, page~2 in \cite{MR1410258}) and the Gaussian mean width of the ``unit ball" associated with $\Psi(\cdot)=\norm{\cdot}_p$ for $0<p<1$ can also be computed.  

To that end, let $\{e_1,\ldots,e_d\}$ be the canonical basis of $\R^d$. Since $\{\pm e_1,\ldots,\pm e_d\}\subset B_p^d\subset B_1^d$ for $p<1$, it is evident that $\ell^*(B_p^d)\sim \sqrt{\log  d}$. Thus, the error rates of the LASSO, obtained in Theorem~\ref{theo:lasso-small-ball-weak-moments}, dominate all the $\ell_p$-regularization rates when $0<p \leq 1$. We therefore obtain the same result for $\ell_p$-regularization with $0<p<1$ as the one in the first case of Theorem~\ref{theo:ell_p_reg} for $\ell_p$-regularization with $1\leq p\leq 1+1/\log p$.

However, the resulting rate is not the minimax rate in the true model, as can be seen from \cite{MR2882274}. Indeed, fix $0<p\leq1$. Consider an unknown $t^*\in \rho B_p^d$ and the corresponding Gaussian linear model $Y_i=\inr{x_i,t^*} + W_i, i=1,\ldots, N$, where the matrix whose rows are $(x_i)_{i=1}^N$ satisfies some RIP property and $W_1, \cdots, W_N$ are independent, centered Gaussian variables with variance $\sigma^2$.  For specific asymptotics of $N$ and $d$ (see \cite{MR2882274} for a precise formulation), the authors show that minimax rate of the problem is given by
\begin{equation*}
\sigma^2 \rho \Big(\frac{\log d}{N}\Big)^{1-\frac{p}{2}},
\end{equation*}
and similar results have been obtained in \cite{WPGaoYang12}. Thus, our estimate recovers the minimax rate in the true model only when $p=1$. When $0<p<1$, it is possible that the choice of the $\Psi(t)=\|t\|_p$ is suboptimal, and instead one should use $\Psi(t)=\norm{t}_p^p$ as was suggested in \cite{MR2816342} for the problem of $S_p$-regularization for $0<p\leq1$.
\end{remark}

\vskip0.4cm

\noindent{\bf Example: weak-$\ell_p$-regularization for $0<p\leq1$.}
Weak-$\ell_p$ norms have been used to model sparsity in High-Dimensional Statistics (see, for instance, \cite{MR2281879, WPGaoYang12}). To define those norms, let $t_1^*\geq t_2^*\geq \ldots \geq t_d^*$ be the
non-increasing rearrangement of $(|t_i|)_{i=1}^d$. Set $\norm{t}_{p\infty}=\max_{1\leq j\leq d} j^{1/p} t^*_j$ and put $B^d_{p\infty}=\{t\in\R^d:t^*_j\leq j^{-1/p} \mbox{ for every } 1\leq
  j\leq d\}$.

One can use the following well-known fact (see, e.g., Theorem~B in \cite{MR2371614}) to control the Gaussian mean-width of the unit ball associated with $\norm{\cdot}_{p\infty}$.
\begin{Proposition}
For $0<p\leq 1$,
\begin{equation*}
 \ell^*(B_{p\infty})\lesssim  \left\{
  \begin{array}[c]{cc}
   \frac{\sqrt{\log d}}{p-1} & \mbox{ if } 0<p<1
   \\
   \\
   \big(\log d\big)^{3/2} & \mbox{ if } p=1.
  \end{array}
\right.
\end{equation*}
\end{Proposition}

Now, one may apply  Theorem~\ref{theo:learning-linear-norm} and obtain the following result.

\begin{Theorem}
Under the assumptions of Theorem \ref{theo:learning-linear-norm} and using its notation,
  \begin{itemize}
  \item  If $p<1$ and
  \begin{equation*}
      \hat t \in\argmin_{t\in
    \R^d}\Big(\frac{1}{N}\sum_{i=1}^N(Y_i-\inr{X_i,t})^2+c_2\eta_p^3\sigma_q L wu\norm{t}_{p\infty}\sqrt{\frac{\log d}{N}}\Big),
  \end{equation*}
  then with probability larger than the one in \eqref{eq:proba1}
  \begin{equation*}
    \norm{\hat t-t^*}_2^2\lesssim_{p,L,q} \left\{
      \begin{array}{cc}
        \sigma_qwu\norm{t^*}_{p\infty}\sqrt{\frac{\log d}{N}} &    \mbox{ if } N
        \gtrsim_L d,
        \\
        \\
        \max\left\{\sigma_qwu\norm{t^*}_p\sqrt{\frac{\log d}{N}},
        \norm{t^*}_{p\infty}^2\frac{\log d}{N}\right\} & \mbox{ otherwise}.
      \end{array}
\right.
  \end{equation*}
\item If $p=1$ and
  \begin{equation*}
      \hat t \in\argmin_{t\in
    \R^d}\Big(\frac{1}{N}\sum_{i=1}^N(Y_i-\inr{X_i,t})^2+
    c_2\eta_1^3\sigma_q\norm{t}_{1\infty}\sqrt{\frac{\log^3 d}{N}}\Big),
  \end{equation*}
then with probability larger than the one in \eqref{eq:proba1}
  \begin{equation*}
    \norm{\hat t-t^*}_2^2\lesssim_{L,\delta,q} \left\{
      \begin{array}{cc}
        \sigma_qwu\norm{t^*}_{1\infty}\sqrt{\frac{\log^3 d}{N}} &    \mbox{ if } N
        \gtrsim_L d,
        \\
        \\
        \max\left\{\sigma_q wu\norm{t^*}_{1\infty}\sqrt{\frac{(\log d)^3}{N}},
        \norm{t^*}_{1\infty}^2\frac{\log^3 d}{N}\right\} & \mbox{ otherwise}.
      \end{array}
\right.
  \end{equation*}
  \end{itemize}
\end{Theorem}

\vskip0.4cm
\noindent{\bf Example: the Micchelli, Morales and Pontil's regularization functions.}

Let $\Theta$ be a nonempty convex cone in $[0,\infty)^d$, and for every $t\in\R^d$ set
\begin{equation}
  \label{eq:MMP-norm} \Omega(t|\Theta)=\inf_{\theta\in\Theta}\frac12\sum_{j=1}^d\Big(\frac{t_j^2}{\theta_j}+\theta_j\Big).
\end{equation}
It was shown in \cite{MMP10} that $\Omega(t|\Theta)$ is a norm on $\R^d$.

This family of norms captures several classical regularization functions, by an appropriate choice of the cone $\Theta$. For instance, the $\ell_1^d$-norm is obtained by selecting
$\Theta=[0,\infty)^d$. Also, the {\it group LASSO} introduced in \cite{MR2212574} is  generated by a cone: indeed, if $(G_1,\cdots,G_T)$ is a partition of $\{1,\ldots,d\}$ and
\begin{equation}\label{eq:theta_group_lasso}
 \Theta=\{\theta\in[0,\infty)^d \mbox{ that is constant
  within each group } G_\ell\},
\end{equation}
then
\begin{equation*}
  \Omega(t|\Theta)=\sum_{\ell=1}^T \sqrt{|G_\ell|}\norm{t_{|G_\ell}}_2,
\end{equation*}
where $|G_\ell|$ is the cardinality of the set of coordinates $G_\ell$ and $t_{|G_\ell}$ is
the restriction of $t$ to $G_\ell$.

Error bounds for procedures that use $\Psi(t)=\Omega(t|\Theta)$ as regularization functions have been established in \cite{MR2913714}, under the assumption that the loss functions is bounded and Lipschitz (see Theorem~1 there). Naturally, the squared loss is not covered by such a result because it is not bounded in $\R^d$, nor is it Lipschitz. Our aim is to provide similar results as the one in \cite{MR2913714} for a quadratic loss for a subgaussian random vector $X$ and a noise in $L_q$ for some $q>2$. To that end, we first compute the Gaussian mean width of the unit ball of such a norm.

\begin{Proposition}
  \label{prop:MMP-gauss-width}
Let $\Theta$ be a nonempty convex cone in $[0,\infty)^d$ and set ${\cal B}=\{t : \Omega(t|\Theta) \leq 1\}$. Let $S_1^{d-1}$ be the unit sphere of $\ell_1^d$ and put ${\cE x}$  to be the set of extreme points of $\Theta\cap S_1^{d-1}$. If $M=\max_{a\in {\cE x}}\norm{a}_\infty^{1/2}$, then, for an absolute constant $c$,
\begin{equation} \label{eq:Theta-ball}
  \ell^*({\cal B})\leq 1+ cM\sqrt{2\log\big(|{\cE x}|\big)}.
\end{equation}
\end{Proposition}

The proof of Proposition~\ref{prop:MMP-gauss-width} may be derived in various ways (see a similar result in \cite{MR2913714}), though we will use a chaining argument which actually leads to a stronger estimate than \eqref{eq:Theta-ball}.

\begin{Definition}
Let $T\subset\R^d$ and $\norm{\cdot}$  be a norm on $\R^d$. For every $\alpha>1$ set
\begin{equation*}
\gamma_\alpha(T,\norm{\cdot})=\inf_{(T_s)}\sup_{t\in T}\sum_{s=0}^\infty 2^{s/\alpha}\norm{\pi_{s+1}t-\pi_st}
\end{equation*}
where the infimum is taken with respect to all sequences $(T_s)$ of subsets of $T$ for which $|T_0|=1$ and for $s \geq 1$, $|T_s|\leq 2^{2^s}$, and $\pi_s t$ is the nearest point to $t$ in $T_s$ with respect to  $\norm{\cdot}$.
\end{Definition}
Clearly, if $T$ is finite then
$\gamma_\alpha(T,\| \ \|) \lesssim \sup_{t \in T} \|t\| \cdot \log^{1/\alpha} |T|.$

\vspace{0.6cm}

\noindent{\bf Proof of Proposition \ref{prop:MMP-gauss-width}.}
It is straightforward to verify (see, e.g., \cite{MMP10}) that the dual norm to $\Omega(\cdot|\Theta)$ is
\begin{equation}
  \label{eq:dual-norm-MMP}
  \Omega^*(t|\Theta)=\max_{a\in{\cE x} }\Big(\sum_{j=1}^d a_j t_j^2\Big)^{1/2}.
\end{equation}
Let $g_1,...,g_d$ be independent, standard Gaussian random variables,
Applying a Bernstein type inequality for a sum of independent $\psi_1$ random variables (see Corollary 2.10 in \cite{MR1269606}), it follows that for every $a_1,...,a_N$, every $u>0$ and any $s\in\N$,
$$
Pr\left[ \Big|\sum_{j=1}^d a_j (g_j^2-1) \Big| \geq u2^{s/2}\|a\|_2 + u^22^s \|a\|_\infty\right] \leq 2\exp(c_12^su^2).
$$
Hence, using a standard chaining argument,
$$
\E \sup_{a \in {\cE x}} \sum_{j=1}^d a_j g_j^2 \leq 1 + c_2 \left(\gamma_2({\cE x}, \| \ \|_2) + \gamma_1({\cE x}, \| \ \|_\infty)\right).
$$
Now one may apply the trivial estimates on $\gamma_1$ and $\gamma_2$. Firstly, $\gamma_1(\cE x, \norm{\cdot}_\infty)\lesssim M^2 \log (|\cE x|)$, and secondly, noting that $|\sum_{j=1}^d a_j| \leq \|a\|_1=1$ and thus $\norm{a}_2\leq \norm{a}_\infty^{1/2}$, one has $\gamma_2(\cE x, \norm{\cdot}_2)\lesssim M \sqrt{\log (|\cE x|)}$. Therefore, by Jensen's inequality,
$$
\E \sup_{a \in \cE x} \bigl(\sum_{j=1}^d a_j g_j^2\bigr)^{1/2} \leq 1 + cM \sqrt{\log (|\cE x|)}.
$$
\endproof

\begin{Theorem}\label{theo:pontil}
Using the notation above and of Theorem \ref{theo:learning-linear-norm}, let
$$
\Lambda(t)=\Omega(t|\Theta) M \sqrt{\frac{\log(|\cE x|)}{N}}.
$$
If
\begin{equation*}
      \hat t \in\argmin_{t\in
    \R^d}\Big(\frac{1}{N}\sum_{i=1}^N(Y_i-\inr{X_i,t})^2+
    c_2\sigma_q L wu \Lambda(t)\Big)
  \end{equation*}
then with probability larger than the one in \eqref{eq:proba1}
  \begin{equation*}
\|\inr{\hat{t}-t^*,\cdot}\|_{L_2(\mu)}^2  \lesssim_{L,q} \left\{
      \begin{array}{cc}
        \sigma_q wu \Lambda(t^*) &    \mbox{ if } N
        \gtrsim_L d,
        \\
        \\
        \max\left\{\sigma_q L wu\Lambda(t^*),
        \Lambda^2(t^*)\right\} & \mbox{ otherwise}.
      \end{array}
\right.
  \end{equation*}
\end{Theorem}

When $\Theta=[0,\infty)^d$ then $M \sqrt{\log({\cE x})} \lesssim \sqrt{\log d}$. Hence, Theorem~\ref{theo:pontil} yields the same error rate as the one obtained for the LASSO in Theorem~\ref{theo:lasso-small-ball-weak-moments} and Theorem~\ref{theo:ell_p_reg}, though under a stronger subgaussian assumption on $X$. This is not surprising because when $\Theta=[0,\infty)^d$, $\Omega(t|\Theta)=\|t\|_1$ and the resulting RERM is just the LASSO.

In the case of the group LASSO, for $\Theta$ as in \eqref{eq:theta_group_lasso}, one has $M \sqrt{\log( {\cE x})} \lesssim \sqrt{\log |T|}$, and $\Lambda(t^*) \sim \Omega(t^*|\Theta) M \sqrt{(\log |T|)/N}$.

\vskip0.5cm
\noindent{\bf Example: The SLOPE regularization}

In \cite{slope2,slope1}, the authors introduced the
regularization function:
\begin{equation*}
\Psi(t)=  \norm{t}_{SLOPE}=\sum_{j=1}^d\lambda_jt_j^*
\end{equation*}where $\lambda_1\geq \cdots\geq \lambda_d\geq0$ and
$t_1^*\geq \cdots\geq t^*_d\geq0$ is the non-increasing rearrangement of $(|t_i|)_{i=1}^d$.

In \cite{slope2} the given data is generated by the Gaussian linear model $Y_i=\inr{X_i, t^*} + W_i, i=1,\ldots,N$ for a Gaussian design $X_i\sim\cN(0, N^{-1}I_{d\times d})$ (note that the covariance matrix is normalized by $1/N$) and a centered Gaussian noise $W_i$ with variance $\sigma^2$ that is independent of the design $X_i$. Setting $\Phi^{-1}(\alpha)$ to be the $\alpha$-th quantile of a standard normal distribution and $q\in(0,1)$, the weights were chosen to be
\begin{equation}
  \label{eq:14}
  \lambda_i=\Phi^{-1}(1-iq/(2d)),
\end{equation}
and, for this choice of weights, SLOPE was defined by
\begin{equation*}
  \hat t \in \argmin_{t\in \R^d}\Big(\frac{1}{2N}\sum_{i=1}^N(Y_i-\inr{X_i,t})^2+\sigma\frac{\norm{t}_{SLOPE}}{\sqrt{N}}\Big).
\end{equation*}

The result in \cite{slope2} is asymptotic in the sample size $N$ and in the dimension $d$ in the following sense:

\begin{Theorem}[Theorem~1.2 \cite{slope2}]\label{theo:candes_slope}
Let $0<\eps<1$ and set $1\leq k\leq d$ that satisfy $k/d=o(1)$ and $(k \log d)/N=o(1)$ when
$N\rightarrow \infty$. Then,
\begin{equation*}
\lim_{N \to \infty}  \sup_{\norm{t^*}_0\leq k}Pr\Big(\frac{N\norm{\hat
      t-t^*}_2^2}{2\sigma^2 k \log(d/k)}>1+3\eps\Big) = 0,
\end{equation*}
where the supremum is taken with respect to all vectors that are supported on at most $k$ coordinates.
\end{Theorem}

It was shown in \cite{slope2} that $2\sigma^2 k \log(d/k)/N$ is the (asymptotic) minimax rate for $t^*$ that is $k$-sparse.

\vskip0.4cm
The article \cite{slope2} (see Section~6 there) raises the question of extending Theorem \ref{theo:candes_slope} beyond the Gaussian case, especially when the coordinates of $X$ may be correlated. We study this question in the context of sparse recovery and for an arbitrary choice of weights in \cite{LM_sparsity}, leading to error bounds that depend on $\|t^*\|_{0}$. Here, we obtain a complexity-dependent error rate that depends on $\|t^*\|_{SLOPE}$.

\begin{Proposition}
Set ${\cal B}=\{t\in\R^d: \norm{t}_{SLOPE}\leq1\}$.
There exists an absolute constant $C$, for which, if $M = \max_{1 \leq j \leq d} \lambda_j^{-1} \sqrt{\log(ed/j)}$, then $\ell^*({\cal B}) \leq C M$.
\end{Proposition}
\proof
The proof is outcome of a standard binomial estimate. Let $G=(g_i)_{i=1}^d$ be a standard Gaussian vector and observe that
\begin{align*}
  &\ell^*({\cal B})=\E\sup_{t\in {\cal B}}\inr{G,t} \leq \E\sup_{t\in {\cal B}}\sum_{j=1}^dg_j^* t_j^*
    \leq \E\sup_{t\in {\cal B}}\sum_{j=1}^d \frac{g_j^*}{\lambda_j} \lambda_jt_j^*\leq \E \max_{1\leq j\leq d}\frac{g_j^*}{\lambda_j}.
\end{align*}
For $u \geq 1$,
\begin{align*}
  &Pr\bigl(\max_{1\leq j\leq d}\frac{g_j^*}{\lambda_j}\geq u\bigr)\leq
  \sum_{j=1}^d Pr\left(g_j^*\geq u \lambda_j\right) \leq \sum_{j=1}^d \binom{d}{j} Pr^j\left(|g|\geq u \lambda_j\right)
  \\
& \leq 2\sum_{j=1}^d
\exp\Big(j\log\Big(\frac{ed}{j}\Big)-c_1ju^2\lambda_j^2\Big) \leq 2\exp(c_2u^2),
\end{align*}
where the last inequality follows if one sets $u^2 \geq \max_j \lambda_j^{-2}\log(ed/j)$. The proof is concluded by integrating the tails.
\endproof

Theorem~\ref{theo:learning-linear-norm} leads to estimation properties
of SLOPE.

\begin{Theorem}\label{theo:slope}
Using the notation of Theorem \ref{theo:learning-linear-norm}, if $\Psi(t)=\|t\|_{SLOPE}$, $\max_j \lambda_j^{-1} \sqrt{\log(ed/j)} \leq C$ and
  \begin{equation*}
      \hat t \in\argmin_{t\in
    \R^d}\Big(\frac{1}{N}\sum_{i=1}^N(Y_i-\inr{X_i,t})^2+c_2\sigma_q Lwu\frac{\norm{t}_{SLOPE}}{\sqrt{N}}\Big),
  \end{equation*}
then with probability larger than the one in \eqref{eq:proba1},
  \begin{equation*}
    \norm{\hat t-t^*}_2^2\lesssim_{L,q,C} \left\{
      \begin{array}{cc}
       \frac{ \sigma_q wu\norm{t^*}_{SLOPE}}{\sqrt{N}} &    \mbox{ if } N
        \gtrsim_L d,
        \\
        \\
        \max\left\{\frac{\sigma_q L wu\norm{t^*}_{SLOPE}}{\sqrt{N}},
        \frac{\norm{t^*}_{SLOPE}^2}{N}\right\} & \mbox{ otherwise}.
      \end{array}
\right.
  \end{equation*}
\end{Theorem}

As was done for the LASSO, one may combine the sparsity-dependent error rate for SLOPE from \cite{LM_sparsity} and the complexity-dependent error rate from Theorem~\ref{theo:slope}. To that end, assume that $X$ is isotropic, $L$-subgaussian and that the noise $\xi$ is in $L_q$ for some $q>2$. Then, with probability larger than the one in \eqref{eq:proba1}
\begin{equation*}
\norm{\hat t-t^*}_2^2\lesssim_{L,q,C} \min\left\{\frac{\sigma_q\norm{t^*}_0}{N}\log\Big(\frac{ed}{\norm{t^*}_0}\Big), \max\left\{\frac{\sigma_q wu\norm{t^*}_{SLOPE}}{\sqrt{N}},
        \frac{\norm{t^*}_{SLOPE}^2}{N}\right\}\right\},
\end{equation*}
for $N\gtrsim \norm{t^*}_0 \log(ed/\norm{t^*}_0)$.

\subsection{Regularization methods in  $\R^{m\times T}$}
\label{sec:regul-meth-matrices}
In this section, we  assume that $X$ takes values in the set of $m\times T$ matrices, endowed with the inner product
$ \inr{A,B}=\sum_{u,v}A_{uv}B_{uv}$. We consider
$ A^*\in\argmin_{A\in\R^{m\times T}}\E\big(Y-\inr{X,A}\big)^2$ and thus $\inr{X,A^*}$ is the best (linear) approximation of $Y$ in the $L_2$ sense.

Let $\Lambda(\rho)\geq \rho \ell^*(K)/\sqrt{N}$ for all $\rho>0$, $u,w>C$ and set
\begin{equation}\label{eq:reg-matrice}
  \hat A \in\argmin_{A\in \R^{m\times     T}}\Big(\frac{1}{N}\sum_{i=1}^N(Y_i-\inr{X_i,A})^2
  +c_2\eta^3\sigma_q L wu \Lambda(\Psi(A))\Big).
\end{equation}
By Theorem~\ref{theo:learning-linear-norm}, with probability larger than the one in \eqref{eq:proba1}
\begin{equation*}
 \norm{\hat A-A}_2^2= \norm{\inr{X,\hat A-A^*}}_{L_2}^2 \lesssim r(10\eta \Psi(A^*))^2
\end{equation*}
where for $\rho\geq0$,
\begin{equation*}
r(\rho)^2\sim_{L,q} \left\{
\begin{array}{cc}
  \sigma_q wu \Lambda(\rho)
 & \mbox{ when } N\gtrsim_L mT
\\
\\
 \max\Big\{\sigma_q wu \Lambda(\rho), \Lambda^2(\rho)\Big\} & \mbox{ otherwise.}
  \end{array}
\right.
\end{equation*}

Let us turn to estimates on $\ell^*(K)$ for the unit balls of the regularization functions used in the {\it matrix completion} and  {\it collaborative filtering} problems.

\vskip0.4cm
\noindent{\bf Example: $S_p$-regularization for $p\geq1$.}

For any $A\in\R^{m\times T}$, let $\sigma_1(A)\geq \sigma_2(A)\geq\cdots\geq \sigma_{m\wedge T}(A)$ be the ordered singular values of $A$ and set $m\wedge T = \min\{m,T\}$. Recall that the $p$-Schatten norm $\norm{\cdot}_{S_p}$ of $A$ is defined by
\begin{equation*}
  \norm{A}_{S_p}=\Big(\sum_{j=1}^{m\wedge T} \sigma_j(A)^p\Big)^{1/p}.
\end{equation*}

Schatten norms have been used extensively in {\it matrix completion} and in {\it collaborative filtering}. Exact reconstruction properties of various procedures have been established via the minimization of the $S_1$-norm, constrained to matching the data (see, e.g.,
\cite{MR2809094,MR2723472,MR2565240,MR2815834,MR2989474}). $S_1$ regularization has also been used in the noisy setup for independent subgaussian noise and, in most case, for subgaussian or deterministic designs, in
\cite{MR2906869,MR2816342,MR2829871,MR2930649,MR2882272,MR3160583}.

\vskip0.4cm

A result that is closely related to ours is Theorem~9.2 from \cite{MR2829871}, in which $X$ is isotropic and $L$-subgaussian; $\xi$ is a symmetric random variable that is independent of $X$ and for which $\|\xi\|_{\psi_\alpha} < \infty$ for some $\alpha \geq 1$ (cf. \cite{MR1113700} for more details on the $\psi_\alpha$-norms); and the target is  $Y=\inr{X,A^*}+\xi$.

Let $N\gtrsim m \cdot {\rm rank}(A^*)$ and set
\begin{equation*}
  \lambda\gtrsim \max\left\{\norm{\xi}_2\sqrt{\frac{m(t+\log m)}{N}},  \norm{\xi}_{\psi_\alpha}\log^{1/\alpha}\Big(\frac{\norm{\xi}_{\psi_\alpha}}{\norm{\xi}_{L_2}}\Big)\frac{\sqrt{m}(t
    + \log N)(t+\log m)}{N}\right\}.
\end{equation*}
The $S_1$-regularization procedure with regularization parameter $\lambda$ satisfies that for every $t>0$, with probability larger than
$1-3\exp(-t)-\exp(-c_0 N)$,
\begin{equation}
  \label{eq:vlad}
  \norm{\hat A-A^*}_{S_2}^2\lesssim \min\left\{\lambda
\norm{A^*}_{S_1}, \lambda^2{\rm rank }(A^*)\right\}.
\end{equation}

In comparison, an estimation result for $S_p$-norm regularization (for any $p\geq1$) follows from Theorem~\ref{theo:learning-linear-norm}, and does not require any assumptions on the ``noise" $\xi=Y-\inr{A^*,X}$, other than $\xi \in L_q$ for some $q>2$. in particular, $\xi$ need not belong to $\psi_\alpha$, nor does it have to be independent of $X$. The result uses the following estimate on the Gaussian mean-width of the unit ball of $S_p$-norms (see, for instance, Proposition~1.4.4 in  \cite{MR3113826}):

\begin{Proposition}
  \label{prop:gaussian-width-schatten}
Let $p\geq1$ and set $B_p^{mT}$ to be the unit ball of
$\norm{\cdot}_{S_p}$. Then
\begin{equation*}
  \ell^*(B_p^{mT})\sim \min\{m,T\}^{1-1/p}\sqrt{m+T}.
\end{equation*}
\end{Proposition}

Combining the previous result with Theorem~\ref{theo:learning-linear-norm}, one obtains the following:
\begin{Theorem}\label{theo:Sp}
Assume that the assumptions  of Theorem~\ref{theo:learning-linear-norm} hold. Let $\Lambda_p(\rho)=\rho \min\{m,T\}^{1-1/p} \sqrt{\frac{m+T}{N}}$ for all $\rho>0$ and
  \begin{equation*}
      \hat A \in\argmin_{A\in
    \R^{m\times T}}\Big(\frac{1}{N}\sum_{i=1}^N(Y_i-\inr{X_i,A})^2+c_2
    \sigma_q L wu\Lambda_p(\norm{A}_{S_p})\Big).
  \end{equation*}
Then with probability at least $1-\delta$,
  \begin{equation*}
 \norm{\hat A-A^*}_{S_2}^2 \lesssim_{p,L,q} \left\{
      \begin{array}{cc}
       \sigma_q wu \Lambda_p(\norm{A^*}_{S_p}) &    \mbox{ if } N
        \gtrsim_L mT,
        \\
        \\
        \max\Big\{\sigma_q wu \Lambda_p(\norm{A^*}_{S_p}),
        \Lambda_p^2(\norm{A^*}_{S_p})\Big\} & \mbox{ otherwise}.
      \end{array}
\right.
  \end{equation*}
\end{Theorem}

\begin{remark}
As in the vector case mentioned earlier, Theorem~\ref{theo:learning-linear-norm} also applies for  $S_p$-regularization for $0<p<1$. In that case, Assumption~\ref{assum:reg-function} is satisfied for $\eta=2^{1/p}$ and the Gaussian mean width of the $S_p$-unit ball satisfies $\ell^*(B_p^{mT})\lesssim \sqrt{m+T}$. It therefore follows from Theorem~\ref{theo:learning-linear-norm} that under the same assumptions as in Theorem~\ref{theo:learning-linear-norm} and for $\Lambda(\rho) = \rho \sqrt{(m+T)/N}$ for all $\rho>0$, the RERM
  \begin{equation*}
      \hat A \in\argmin_{A\in
    \R^{m\times T}}\Big(\frac{1}{N}\sum_{i=1}^N(Y_i-\inr{X_i,A})^2+c_2
  \sigma_q L wu \Lambda(\norm{A}_{S_p})\Big),
  \end{equation*}
satisfies, with probability larger than the one in \eqref{eq:proba1},
  \begin{equation*}
    \norm{\hat A-A^*}_{S_2}^2\lesssim_{p,L,q} \left\{
      \begin{array}{cc}
       \sigma_q wu \Lambda(\norm{A^*}_{S_p}) &    \mbox{ if } N
        \gtrsim_L mT,
        \\
        \\
        \max\Big\{\sigma_qwu\Lambda(\norm{A^*}_{S_p}),
        \Lambda^2(\norm{A^*}_{S_p})\Big\} & \mbox{ otherwise}.
      \end{array}
\right.
  \end{equation*}

Observe that just as in the vector case, when $0<p<1$ this rate is not the minimax rate in the true model  $\norm{A^*}_{S_p}B_p^{mT}$. Indeed, \cite{MR2816342} provides the minimax rate, and, in fact, also shows that the minimax rate may be attained using $\Psi(A)=\norm{A}_{S_p}^p$ as a regularization function. To be more accurate, \cite{MR2816342} considers the following problem: let $x_1,\ldots, x_N$ be $N$ deterministic matrices in $\R^{m\times T}$ satisfying some RIP property and set $W_i$ to be $N$ independent, standard Gaussian variables with variance $\sigma^2$. Set $Y_i=\inr{x_i,A^*} + W_i, i=1,\ldots,N$, leading to the so-called matrix regression model with Gaussian noise and a deterministic design. It is shown in \cite{MR2816342} that when $\rho B_p^{mT}$ for some $0<p\leq1$, the minimax rate of the problem in $\rho B_p^{mT}$ is
\begin{equation*}
\sigma^2 \rho^p\Big(\frac{m+T}{N}\Big)^{1-\frac{p}{2}}
\end{equation*}
in some specific range of $N,\sigma$ and $\rho$. Our result recovers this rate only for $p=1$. 
\end{remark}

\vskip0.4cm
\noindent{\bf Example: Max-norm regularization.}
The max-norm of a matrix is defined by
\begin{equation*}
\norm{A}_{max}=\inf_{A=U V^\top}\norm{U}_{2\rightarrow \infty}\norm{V}_{2\rightarrow \infty},
\end{equation*}
with the infimum is taken with respect to all pairs of matrices $U,V$ for which $A=U V^\top$.

Constrained empirical risk minimization procedures that are based on the max-norm have been used in \cite{CaiZhou13,SS05} for bounded and Lipschitz loss functions and in \cite{LM13} for the squared loss and for a subgaussian and isotropic design vector $X$ and a subgaussian noise $\xi$ independent of $X$. One may show that the minimax rate in the matrix regression model $Y_i=\inr{X_i, A^*} + W_i, i=1, \ldots, N$ where $X_1,\ldots, X_N$ are independent isotropic and subgaussian matrices, $W_1, \ldots, W_N$ are independent centered gaussian variables with variance $\sigma^2$ that are independent of the $X_i$'s and $A^*$ belongs to  the max-norm ball of radius $\rho$, is
\begin{equation}\label{eq:minimax_max_norm}
\max\left\{\sigma \rho \sqrt{\frac{(mT)(m+T)}{N}}, \frac{\rho^2(mt)(m+T)}{N}\right\}
\end{equation} for some specific regime of $\rho, \sigma$ and $N$ (cf. \cite{LM13}).

 To apply Theorem~\ref{theo:learning-linear-norm}, let us estimate the Gaussian mean-width of the unit ball of the max-norm ball, that is, of ${\cal B}=\{A\in\R^{m\times T}: \norm{A}_{max}\leq1\}$.
\begin{Lemma}
There exists an absolute constant $c$ for which, for every $m$ and $T$,
\begin{equation*}
  \ell^*({\cal B})\lesssim \sqrt{(mT)(m+T)}.
\end{equation*}
\end{Lemma}

\proof An application of Grothendieck's inequality (see, e.g., \cite{SS05}) shows that
\begin{equation*}
  {\rm conv}\big(\cX_{\pm}\big)\subset {\cal B}\subset K_G {\rm conv}\big(\cX_{\pm}\big)
\end{equation*}
where $K_G$ is the Grothendieck constant and
$\cX_{\pm}=\{uv^\top:u\in\{\pm1\}^m, v\in\{\pm1\}^T\}$. If
$\mathfrak{G}=(g_{ij})_{1\leq u\leq m:1\leq v\leq T}$ is a standard  $m\times T$
Gaussian matrix, it follows from the Gaussian maximal
inequality (see, e.g., Chapter~3 in \cite{LT:91}) that
\begin{align*}
  &\ell^*({\cal B}) = \E\sup_{A \in {\cal B}}|\inr{\mathfrak{G},A}| \leq
  K_G \E\sup_{A\in {\rm conv}(\cX_{\pm})}|\inr{\mathfrak{G},A}|\\
&= K_G \E\sup_{A\in \cX_{\pm}}|\inr{\mathfrak{G},A}| \lesssim \max_{A\in\cX_{\pm}}\norm{A}_{HS}\sqrt{\log|\cX_{\pm}|} \lesssim \sqrt{(mT)(m+T)}.
\end{align*}
\endproof

\begin{Theorem}\label{theo:max_norm}
Using the assumptions and notation of Theorem~\ref{theo:learning-linear-norm}, and setting $\Lambda(\rho)=\rho \sqrt{(mT)(m+T)/N}$, if
  \begin{equation*}
      \hat A \in\argmin_{A\in
    \R^{m\times T}}\Big(\frac{1}{N}\sum_{i=1}^N(Y_i-\inr{X_i,A})^2+c_2 \sigma_q Lwu \Lambda(\norm{A}_{max})\Big),
  \end{equation*}
  then with probability larger than the one in \eqref{eq:proba1}
  \begin{equation*}
\|\hat{A}-A^*\|_{S_2}^2   \lesssim_{L,q} \left\{
      \begin{array}{cc}
        \sigma_q wu \Lambda(\norm{A^*}_{max}) &    \mbox{ if } N
        \gtrsim_L mT,
        \\
        \\
        \max\Big\{\sigma_q wu \Lambda(\norm{A^*}_{max}), \Lambda^2(\norm{A^*}_{max})\Big\} & \mbox{ otherwise}.
      \end{array}
\right.
  \end{equation*}
\end{Theorem}

As a consequence, we recover the minimax rate of convergence in the matrix regression model with subgaussian design and gaussian noise in the class $\norm{A^*}_{max}{\cal B}$ thanks to max-norm regularization and without knowing $\norm{A^*}_{max}$ in advance.

\vskip0.5cm
\noindent{\bf Example: Atomic-norm regularization.}

The atomic-norm has been used in \cite{MR2989474} in the context of exact and robust recovery using few Gaussian linear measurements of a signal or of a matrix.

Given $\cA\subset\R^{m\times T}$, the elements in $\cA$ are called \textit{atoms}. Set ${\rm conv}(\cA)$ to be the convex hull of $\cA$ and put
\begin{equation}
  \label{eq:gauge}
  \norm{A}_{\cA}=\inf\left\{t>0:A\in t {\rm conv}(\cA)\right\}.
\end{equation}
Even though $ \norm{\cdot}_{\cA}$ need not be a norm (because ${\rm conv}(\cA)$ need not be centrally-symmetric), it is positive homogeneous and satisfies a triangle inequality: for every $A,B\in \R^{m\times T}$ and $\lambda\geq0$:
\begin{equation*}
  \norm{A+B}_\cA\leq \norm{A}_\cA+\norm{B}_\cA \mbox{ and }
  \norm{\lambda A}_\cA=\lambda \norm{A}_\cA.
\end{equation*}
And, if we assume that ${\cA}$ is centrally-symmetric, then $\| \ \|_{\cA}$ is a norm, \eqref{eq:condition-norm} is satisfied and Theorem~\ref{theo:learning-linear-norm} applies.

Set ${\cal B}$ to be the unit ball with respect to $\| \ \|_{\cA}$ and note that $\ell^*({\cal B})=\ell^*(\cA)$. For example, assume that $m=T$ and put $\cA$ to be the set of all orthogonal matrices. Since the unit ball of the spectral norm is the convex hull of the set of orthogonal matrices, one has $\norm{\cdot}_{\cA}=\norm{\cdot}_{S_2}$ and
$$
\ell^*({\cal B})=\E\norm{\mathfrak{G}}_{S_2}\leq
\sqrt{m}\E\norm{\mathfrak{G}}_{S_\infty}\lesssim m.
$$

\begin{Theorem}\label{theo:atom}
Using the assumptions and notation of Theorem~\ref{theo:learning-linear-norm}, let
$\cA\subset\R^{m\times T}$ be a symmetric set of atoms and set $\Lambda(\rho) \geq  \rho\ell^*(\cA)/\sqrt{N}$ for any $\rho>0$. If
  \begin{equation*}
      \hat A \in\argmin_{A\in
    \R^{m\times T}}\Big(\frac{1}{N}\sum_{i=1}^N(Y_i-\inr{X_i,A})^2+c_2
    \sigma_q Lwu \Lambda(\norm{A}_\cA)\Big)
  \end{equation*}
then with probability larger than the one in \eqref{eq:proba1}
  \begin{equation*}
\norm{\hat A-A^*}_{S_2}^2 \lesssim_{L,q}  \left\{
      \begin{array}{cc}
        \sigma_q wu \Lambda(\norm{A^*}_{\cA}) &    \mbox{ if } N
        \gtrsim_L mT,
        \\
        \\
        \max\Big\{\sigma_q wu \Lambda(\norm{A^*}_{\cA}), \Lambda^2(\norm{A^*}_{\cA})\Big\} & \mbox{ otherwise}.
      \end{array}
\right.
  \end{equation*}
\end{Theorem}

\section{Conclusions} \label{sec:conclusion}
We have presented a general result (Theorem~\ref{thm:main}) describing
statistical properties of a constrained regularized procedure in the
learning theoretical framework. This result highlights the role played
by the quadratic and multiplier processes in calibrating the
regularization parameter $\lambda$ as well as their effect on the estimation error rate.
It appears that:
\begin{enumerate}
\item the rates of convergence depend on $\Psi(f^*)$ and we recover the minimax rate in the ``true model" $\{f\in F: \Psi(f)\leq \Psi(f^*)\}$ -- up to a logarithmic factor -- in many well-studied examples .
\item no statistical model is needed to study RERM; all the analysis has been carried out here in the general learning theory setup, and thus without assuming any statistical model. Theorem~\ref{thm:main} and all its corollaries and applications are true regardless of any relation between the target $Y$ and the input $X$. For instance, when predicting $Y$ using linear functionals of $X$ there is no need to assume that $Y$ equals a linear functional of $X$ plus an independent noise; our results hold even if $Y$ were, for instance, a noisy version of a quadratic function of a linear functional of $X$ (e.g. phase retrieval) or even when $Y$ is independent of $X$. 
    
    Our analysis shows that despite considering the more general learning theory framework, the error rate and the regularization parameter used to construct RERM almost match the ones that would have been obtained with more information -- namely, a given statistical model. In the examples we considered, Statistical models are superfluous for the analysis of RERM and as a consequence, they may actually hide what really determines the error rate and the right choice of a regularization parameter:
\begin{itemize}
  \item calibration of the regularization parameter depends only on the
  multiplier process -- which measures the empirical correlation between
  the noise $Y-f^*(X)$ and the class $F$. When this correlation is
  small or even null (in the free-noise case) the regularization
  parameter will also be small.
  \item the key parameters are the ``structure'' of the ``unit ball" of the regularization function (measured here using the Gaussian mean width) and the ``noise level", which we measure through the $L_q$ norm of  $Y-f^*(X)$.
\end{itemize}
\end{enumerate}

\section{Proof of Theorem \ref{thm:limited-moments}} \label{sec:proof-3.6}
Following \cite{shahar_multi_pro}, the proof of Theorem \ref{thm:limited-moments} is based on properties of the following norm:
\begin{Definition} \label{def:(p)-norm}
For a random variable $Z$ and $p \geq 1$, set
$$
\|Z\|_{(p)} = \sup_{1 \leq q \leq p} \frac{\|Z\|_{L_q}}{\sqrt{q}}.
$$
\end{Definition}
The $\norm{\cdot}_{(p)}$ norm is a `local' version of the $\psi_2$
norm. While
\begin{equation*}
\|Z\|_{\psi_2} \sim \sup_{q \geq 1}
\frac{\|Z\|_{L_q}}{\sqrt{q}},
\end{equation*}
$\|Z\|_{(p)}$ captures the subgaussian
behavior of $Z$ up to the $p$-th moment.

Under Assumption~\ref{ass:max-finite}, a high probability bound on
(\ref{eq:finite-max-1}) can be derived from the next result.

\begin{Proposition}[Lemma~2.8 in \cite{LM_compressed}]\label{prop:shahar-lemma-6.4}
There exists an absolute constant $c_0$ for which the following holds.  Let $Z$ be a mean-zero real-valued  random variable and let $Z_1,...,Z_N$ be independent copies of $Z$. Let $p_1 \geq 1$ and assume that   $\|Z\|_{(p_1)}\leq L$, then
  \begin{equation*}
    \norm{\frac{1}{\sqrt{N}}\sum_{i=1}^N Z_i}_{(p_1)}\leq c_0L.
  \end{equation*}
\end{Proposition}

\vskip0.3cm
Setting $U_j = N^{-1/2}\sum_{i=1}^N \eps_i X_i(j)$ and $p_1=\log d$ (recalling that $t\geq1$ and $d\geq N$ in Assumption~\ref{ass:max-finite}), it follows from Proposition~\ref{prop:shahar-lemma-6.4}  that
$$
\|U_j\|_{L_{p_1}} \leq c_0L\sqrt{p_1} \|x_j\|_{L_2}.
$$
Therefore,
\begin{align*}
  &Pr\left(\max_{1\leq j\leq d}|U_j|\geq u\right) \leq \sum_{j=1}^d Pr\left(|U_j|\geq
  u\right)\leq \sum_{j=1}^d \left(\frac{\norm{U_j}_{L_{p_1}}}{u}\right)^{p_1}
  \\
& \leq d \left(\frac{c_0L \sqrt{p_1}\max_{1\leq j \leq d}\norm{x_j}_{L_2}}{u}\right)^{p_1}=d \left(\frac{c_0L \sqrt{\log d}\max_{1\leq j \leq d}\norm{x_j}_{L_2}}{u}\right)^{\log d}.
\end{align*}
Let $w\geq e $ and set $u = c_0 L w \sqrt{\log d}\max_{1\leq j \leq d}\norm{x_j}_{L_2}$; therefore,
\begin{equation} \label{eq:probab-max-of-j-simple}
  Pr\left(\max_{1\leq j\leq d}|U_j|\geq c_1 L w \sqrt{\log d}\max_{1\leq j \leq d}\norm{x_j}_{L_2} \right) \leq
\left(\frac{e}{w}\right)^{\log d},
\end{equation}
which is a high probability estimate on \eqref{eq:finite-max-1} under a limited moment assumption. Integrating the tail,
\begin{equation*}
\E \max_{1\leq j\leq d} \left|\frac{1}{\sqrt{N}}\sum_{i=1}^N \eps_i X_i(j) \right| \lesssim L \sqrt{\log d}\max_{1\leq j \leq d}\norm{x_j}_{L_2}
\end{equation*}
proving \eqref{eq:expected_quad}.

Next, we obtain high probability bounds on \eqref{eq:finite-max-2} -- which requires some preparation.

Let $j\in\{1,\ldots,d\}$ and set $Z_i = X_i(j)$. Consider the Bernoulli sums
$$
Q_j = \sum_{i=1}^N \eps_i \xi_i X_{i}(j)=\sum_{i=1}^N \eps_i \xi_i Z_i.
$$
Denote by $(a_i^*)_{i=1}^N$ the non-increasing rearrangement of $(|a_i|)_{i=1}^N$. A straightforward application of H\"{o}ffding's inequality shows that conditioned on $(\xi_i)_{i=1}^N$ and $(Z_i)_{i=1}^N$, for any $v>0$, with probability at least  $1-2\exp(-v^2/2)$ relative to $(\eps_i)_{i=1}^N$,
\begin{align}\label{eq:hoeffding}
&|Q_j| \leq  \sum_{i \leq m} \xi_i^* Z^*_i + v\left(\sum_{i \geq m} (\xi_i^* Z_i^*)^2\right)^{1/2}
\\
\nonumber \leq & \left(\sum_{i \leq m} (\xi_i^*)^2\right)^{1/2} \left(\sum_{i \leq m} (Z_i^*)^2\right)^{1/2} + v \left(\sum_{i \geq m} (\xi_i^*)^{2r}\right)^{1/2r}\left(\sum_{i \geq m} (Z_i^*)^{2r^\prime}\right)^{1/2r^\prime},
\end{align}
where $r$ and $r^\prime$ are conjugate indices.

As a consequence, high probability bounds on the rearrangements $(\xi_i^*)$ and $(Z_i^*)$ can be used to obtain high probability bounds on $|Q_j|$ (and therefore, on $\max_{1\leq j\leq d}|Q_j|$ as well, using the union bound).


The next two observations, whose proofs may be found in \cite{shahar_multi_pro} give information on the structure of a typical $(Z_i)_{i=1}^N$ when $Z$ has at least $t \log d$ subgaussian moments. It turns out that one may decompose $(Z_i)_{i=1}^N$ to a sum of two vectors, supported on disjoint sets: one consists of the largest $m$ coordinates of $(|Z_i|)_{i=1}^N$, and its $\ell_2^N$ norm is determined by relatively high moments of $Z$; the other one consists of the $N-m$ smaller coordinates of $(|Z_i|)_{i=1}^N$, and if $Z \in L_{q_1}$, its $\ell_r^N$ norm is well-behaved for $r<q_1$.

The `level' $m$ depends on the desired probability estimate and on the moments of $Z$: if one wishes to obtain a probability estimate of $1-2\exp(-p)$ for $p \geq \log N$ (as we will), then $Z$ should have roughly $p$ moments and one should select $m \sim p/\log(eN/p)$.

First, let us consider the smaller coordinates:

\begin{Lemma}[Lemma~3.2 in \cite{shahar_multi_pro}]\label{lemma:single-small-coordinates}
There exist absolute constants $a_0$ and $c_1$ for which the following holds. Let $1 \leq r_1 < q_1$, set $Z \in L_{q_1}$ and put $Z_1,...,Z_N$ to be independent copies of $Z$. Fix $1 \leq p \leq N$, let $u>2$ and set
$$
m = \left\lceil \frac{a_0p}{((q_1/r_1)-1) \log(4+eN/p)} \right\rceil.
$$
If $m>1$, then, with probability at least $1-2u^{-mq_1}\exp(-p)$,
$$
\left(\sum_{i=m}^{N} (Z_i^*)^{r_1} \right)^{1/r_1} \leq c_1\left(\frac{q_1}{q_1-r_1}\right)^{1/r_1} u N^{1/r_1} \|Z\|_{L_{q_1}}
$$and, if $m=1$ and $0<\beta<q_1/r_1-1$ then with probability at least $1-c_2 u^{-q_1} N^{-\beta}$,
\begin{equation*}
\left(\sum_{i=1}^N |Z_i|^{r_1}\right)^{1/r_1}\leq c_1 \Big(\frac{q_1}{q_1-(\beta+1)r_1}\Big) u \norm{Z}_{L_{q_1}} N^{1/r_1}.
\end{equation*}
\end{Lemma}

Next, we consider the larger coordinates:

\begin{Lemma}[Lemma~3.4 in \cite{shahar_multi_pro}]\label{lemma:large-coordinates} There exists absolute constants $a_1$ and $c_0$ for which the following holds. Let $Z_1,...,Z_N$ be independent copies of a random variable $Z$, set $p \geq \log N$ and put $1 \leq m \leq N/2e$ that satisfy $m \leq a_1 p/\log\big(eN/p\big)$. Then,  for every $u>1$, with probability at least $1-u^{-2p}\exp(-p)$, one has
$$
\left(\sum_{i = 1}^m (Z_i^*)^2\right)^{1/2} \leq c_0 u\sqrt{p}\|Z\|_{(2p)}.
$$
\end{Lemma}

In particular, under Assumption~\ref{ass:max-finite}, we apply Lemma~\ref{lemma:single-small-coordinates} and Lemma~\ref{lemma:large-coordinates} to $p,q_1$ and  $r_1$ defined by
\begin{equation*}
2p = t \log d,r_1=2r^\prime \mbox{ and } q_1=r_1\max\Big\{2, 1+\frac{a_0}{a_1}\Big\},
\end{equation*}
where $a_0$ and $a_1$ are the absolute constants from Lemma~\ref{lemma:single-small-coordinates} and \ref{lemma:large-coordinates}. We also set
\begin{equation}\label{eq:def_m}
m= \left\lceil \frac{a_0p}{\log(4+eN/p)} \right\rceil
\end{equation} and observe that if
\begin{equation*}
m_0 = \left\lceil \frac{a_0p}{((q_1/r_1)-1) \log(4+eN/p)} \right\rceil \mbox{ then } m_0\leq m \leq \frac{a_1 p}{\log(eN/p)}
\end{equation*}
and $m_0 q_1 \sim m$. Moreover, if $\kappa_0$ is a large enough absolute constant and  $t\geq \kappa_0$, then $m_0>1$. Recalling that $p \geq 2\log d$ and setting $Z_i=X_i(j)$ for $i=1,\ldots,N$, it follows that for any $u>2$, with probability at least $1-u^{-2p}\exp(-p/2)-2 u^{-c_0m}\exp(-p/2)$, for every $1 \leq j \leq d$
\begin{equation}\label{eq:largest}
\left(\sum_{i=1}^m (Z_i^*)^2\right)^{1/2} \lesssim u \sqrt{p}\norm{Z}_{(2p)}\lesssim uL \sqrt{t \log d} \norm{x_j}_{L_2}
\end{equation} and
\begin{equation}\label{eq:smallest}
\left(\sum_{i= m}^N (Z_i^*)^{2r^\prime}\right)^{1/2r^\prime} \lesssim u \norm{Z}_{L_{q_1}} N^{1/2r^\prime}\lesssim u L \sqrt{r^\prime} \norm{x_j}_{L_2} N^{1/2r^\prime}.
\end{equation}

Let $\xi_1,...,\xi_N$ be independent copies of $\xi$ and recall that $(\xi_i^*)_{i=1}^N$ is the monotone non-increasing rearrangement of $(|\xi_i|)_{i=1}^N$. We apply Lemma~\ref{lemma:single-small-coordinates} for $q_1=q$, $r_1=2r$ and set
\begin{equation*}
 m_1 = \left\lceil\frac{a_0p}{((q_1/r_1)-1) \log(4+eN/p)}\right\rceil.
\end{equation*}
Thus, $m_1>1$ when $t\geq \kappa_0$ for a large enough constant $\kappa_0$, and if $m$ is as in \eqref{eq:def_m} one has  $m\geq m_1$ and $m_1 q_1\sim m$. Hence, for $p=(t/2)\log d$, with probability larger than $1-2 u^{-c_0m}\exp(- p)$,
\begin{equation}\label{eq:small_xi}
\left(\sum_{i= m}^N (\xi^*_i)^{2r}\right)^{1/2r}\leq c(q) u \norm{\xi}_{L_q} N^{1/2r}.
\end{equation}
This provides a high probability bound on the smaller coefficients of $(|\xi_i|)_{i=1}^N$, and now we shall turn to a result on the larger ones.

\begin{Lemma}[Lemma~4.3 in \cite{shahar_multi_pro}]\label{lemma:ell2_xi}
Let $q>2$ and assume that $\xi\in L_q$. If $\xi_1,\ldots,\xi_N$ are independent copies of $\xi$ then for every $w>1$ with probability larger than $1-c_0 w^{-q}N^{-((q/2)-1)}\log^q N$,
\begin{equation*}
\left(\sum_{i=1}^m \xi_i^2\right)^{1/2}\leq \left(\sum_{i=1}^N \xi_i^2\right)^{1/2}\leq c_1 w \norm{\xi}_{L_q}\sqrt{N}.
\end{equation*}
\end{Lemma}


Setting $Z_i=X_i(j)$ for $i=1,\ldots,N$ and applying \eqref{eq:largest}, \eqref{eq:smallest}, \eqref{eq:small_xi} and Lemma~\ref{lemma:ell2_xi}, we obtain, with probability larger than
\begin{equation*}
1-\frac{\exp(-p/2)}{u^{2p}}-\frac{4 \exp(-p/2)}{u^{c_0m}}-\frac{c_0 \log^q N }{w^q N^{q/2-1}}
\end{equation*}
that for every $j=1,\ldots,d$
\begin{equation*}
\left(\sum_{i \leq m} (\xi_i^*)^2\right)^{1/2}\left(\sum_{i \leq m} (Z_i^*)^2\right)^{1/2} \lesssim  uw L \sqrt{t \log d}\norm{x_j}_{L_2} \norm{\xi}_{L_q}\sqrt{N}
\end{equation*}
and
$$
\left(\sum_{i \geq m} (\xi_i^*)^{2r}\right)^{1/2r}\left(\sum_{i \geq m} (Z_i^*)^{2r^\prime}\right)^{1/2r^\prime} \leq c(q) u^2 L \norm{\xi}_{L_q}\sqrt{N}\norm{x_j}_{L_2}.
$$
Then, by plugging those inequalities in \eqref{eq:hoeffding}, it is evident that under Assumption~\ref{ass:max-finite}, for $u>2, v>0$, $w \geq 2$, $2p=t\log d$ and $m\sim p/\log(eN/p)$, with probability at least
\begin{equation*}
1-\frac{\exp(-p/2)}{u^{2p}}-\frac{4 \exp(-p/2)}{u^{c_0 m}}-\frac{c_1 \log^q N }{w^q N^{q/2-1}}-2 \exp(-v^2t\log d),
\end{equation*}
$$
\max_{1\leq j\leq d}\Big|\sum_{i=1}^N \eps_i \xi_i X_i(j)\Big|=\max_{1 \leq j \leq d} |Q_j| \lesssim_q (uw + u^2 v) L\|\xi\|_{L_q}  \sqrt{N}\sqrt{t\log d} \max_{1 \leq j \leq d} \|x_j\|_{L_2}.
$$

\begin{footnotesize}
\bibliographystyle{plain}
\bibliography{biblio}
\end{footnotesize}

\end{document}